\documentclass[reqno,11pt]{amsart}
 
 \usepackage{graphics, color, amsmath, amsfonts, amssymb, amsthm, amscd    
}   
   
\DeclareMathSymbol{\leqslant}{\mathalpha}{AMSa}{"36} 
   
\DeclareMathSymbol{\geqslant}{\mathalpha}{AMSa}{"3E} 
   
\DeclareMathSymbol{\eset}{\mathalpha}{AMSb}{"3F}     
   
\renewcommand{\leq}{\;\leqslant\;}                   
   
\renewcommand{\geq}{\;\geqslant\;}                   
   
\newcommand{\dd}{\,\text{\rm d}}             
   
\DeclareMathOperator*{\union}{\bigcup}       
   
\newcommand{\sumtwo}[2]{\sum_{\substack{#1 \\ #2}}} 
   
\newcommand{\prodtwo}[2]{\prod_{\substack{#1 \\ #2}}}     
   
\makeatletter   
\def\captionfont@{\footnotesize}   
\def\captionheadfont@{\scshape}   
\long\def\@makecaption#1#2{%
  \vspace{2mm}   
  \setbox\@tempboxa\vbox{\color@setgroup   
    \advance\hsize-6pc\noindent   
    \captionfont@\captionheadfont@#1\@xp\@ifnotempty\@xp   
        {\@cdr#2\@nil}{.\captionfont@\upshape\enspace#2}%
    \unskip\kern-6pc\par   
    \global\setbox\@ne\lastbox\color@endgroup}%
  \ifhbox\@ne 
    \setbox\@ne\hbox{\unhbox\@ne\unskip\unskip\unpenalty\unkern}%
  \fi   
  \ifdim\wd\@tempboxa=\z@ 
    \setbox\@ne\hbox to\columnwidth{\hss\kern-6pc\box\@ne\hss}%
  \else 
    \setbox\@ne\vbox{\unvbox\@tempboxa\parskip\z@skip   
        \noindent\unhbox\@ne\advance\hsize-6pc\par}%
\fi   
  \ifnum\@tempcnta<64 
    \addvspace\abovecaptionskip   
    \moveright 3pc\box\@ne   
  \else 
    \moveright 3pc\box\@ne   
\nobreak   
\vskip\belowcaptionskip   
\fi   
\relax   
}   
\makeatother   
   
\def\writefig#1 #2 #3 {\rlap{\kern #1 truecm   
\raise #2 truecm \hbox{#3}}}

\newtheorem{lem}{Lemma}[section]   
   
\newtheorem{pro}{Proposition}[section]   
   
\newtheorem{thm}{Theorem}[section]

\newtheorem{rem}{Remark}[section]

\newcommand{\cA}{\ensuremath{\mathcal A}}

\newcommand{\cC}{\ensuremath{\mathcal C}}

\newcommand{\cE}{\ensuremath{\mathcal E}}   
   
\newcommand{\cF}{\ensuremath{\mathcal F}}   
   
\newcommand{\cG}{\ensuremath{\mathcal G}}   
   
\newcommand{\cH}{\ensuremath{\mathcal H}}

\newcommand{\cN}{\ensuremath{\mathcal N}}

\newcommand{\cS}{\ensuremath{\mathcal S}}

\newcommand{\cU}{\ensuremath{\mathcal U}}   
   
\newcommand{\cV}{\ensuremath{\mathcal V}}

\newcommand{\frn}{\ensuremath{\mathfrak n}}

\newcommand{\frs}{\ensuremath{\mathfrak s}}

\newcommand{\bbC}{{\ensuremath{\mathbb C}} }

\newcommand{\bbE}{{\ensuremath{\mathbb E}} }

\newcommand{\bbH}{{\ensuremath{\mathbb H}} }

\newcommand{\bbN}{{\ensuremath{\mathbb N}} }

\newcommand{\bbP}{{\ensuremath{\mathbb P}} }   
   
\newcommand{\bbQ}{{\ensuremath{\mathbb Q}} }   
   
\newcommand{\bbR}{{\ensuremath{\mathbb R}} }   
   
\newcommand{\bbS}{{\ensuremath{\mathbb S}} }

\newcommand{\bbV}{{\ensuremath{\mathbb V}} }

\newcommand{\bbZ}{{\ensuremath{\mathbb Z}} }

\newcommand{\gb}{\beta}   
   
\newcommand{\gga}{\gamma}            
   
\newcommand{\gd}{\delta}   
   
\newcommand{\gD}{\Delta}

\newcommand{\gl}{\lambda}

\newcommand{\gs}{\sigma}

\newcommand{\ux}{\underline{x}}   
   
\newcommand{\uy}{\underline{y}}   
   
\newcommand{\uz}{\underline{z}}   
   
\newcommand{\ugga}{\underline{\gga}}   
   
\newcommand{\ugl}{\underline{\gl}}

\newcommand{\uempty}{\underline{\emptyset}}

\newcommand{\norma}[1]{\|#1\|}   
\newcommand{\norm}[1]{|#1|}

\newcommand{\normsup}[1]{\|#1\|_{{\scriptscriptstyle\infty}}}   
   
\newcommand{\normgt}[1]{\norma{#1}_{\theta}}   
   
\newcommand{\snormgt}[1]{|\!\!| #1|\!\!|_{\theta}}   
   
\newcommand{\var}[2]{{\rm {\bf var}}_{#1}\left(#2\right)}

\newcommand{\Space}{\cS_\emptyset}   
   
\newcommand{\SpaceN}{\cS_\emptyset^{(N)}}   
   
\newcommand{\Ggt}{\cG_\theta}   
   
\newcommand{\GgtN}{\cG_\theta^{(N)}}   
   
\newcommand{\Qn}[3]{\bbQ_{n, #1}^{#2}\lb #3\rb}   
   
\newcommand{\Qnx}{\bbQ_{n,\ux}}   
   
\newcommand{\Qnxi}[1]{\bbQ_{n,\ux}^{#1}}

\newcommand{\Enx}{\bbE_{n,\ux}}

\newcommand{\Hnx}{\bbH_{n,\ux}}   
   
\newcommand{\chnx}{\widehat{\bbQ}_{n,\ux}}

\def\1{\ifmmode {1\hskip -3pt \rm{I}} \else {\hbox {$1\hskip -3pt \rm{I}$}}\fi} 

\newcommand{\smallo}{o}   
   
\newcommand{\so}{\smallo (1)}   
   
\newcommand{\df}{\stackrel{\Delta}{=}}   
   
\newcommand{\sep}{~|~}

\newcommand{\Kbeta}{{\bf K}_\beta}   
   
\newcommand{\pKbeta}{\partial\Kbeta}   
   
\newcommand{\Ubeta}{{\bf U}_\beta}

\newcommand{\step}[1]{S{\small TEP}~#1}   
   
\newcommand{\case}[1]{C{\small ASE}~#1}   
   
\newcommand{\weight}[1]{q_\beta \lb #1\rb}   
   
\newcommand{\cone}[2]{\cC_{#1}\lb #2 \rb}   
   
\newcommand{\nback}[1]{\#^{{\rm back}}_{#1}}   
   
\newcommand{\nmark}[1]{\#^{{\rm mark}}_{#1}}   
   
\newcommand{\gnz}{g^n_{\tau,\uz}}

\newcommand{\bk}[1]{\langle#1\rangle} 
   
\newcommand{\abs}[1]{\lvert#1\rvert}  
   
\newcommand{\bigabs}[1]{\Bigl\lvert#1\Bigr\rvert} 
   
\newcommand{\ra}{\rightarrow}   
   
\newcommand{\setof}[2]{\left\{#1\,:\,#2\right\}}   
   
\newcommand{\bigsetof}[2]{\Bigl\{#1\,:\,#2\Bigr\}}

\newcommand{\lb}{\left(}   
   
\newcommand{\rb}{\right)}   
   
\newcommand{\lbr}{\left\{}   
   
\newcommand{\rbr}{\right\}}

\newcommand{\symdiff}{{\scriptstyle \triangle}\,}   
   
\newtheorem{bigthm}{Theorem}   
   
\newcommand{\ind}[2]{{\rm ind}(#1,#2)}   
   
\newcommand{\uJ}{{\boldsymbol J}}

\begin{document}   
   
 \title[\,]{Ornstein-Zernike Theory for  
the finite range Ising models above $T_c$}   
\footnotetext{AMS 1991 Subject Classification: 60F15, 60K15, 60K35, 82B20, 37C30\\   
Key Words and Phrases: Ising model, Ornstein-Zernike decay of correlations,   
  Ruelle operator, renormalization, local limit theorems}   
\author{Massimo Campanino}   
\address{   
Dipartimento di Matematica\\   
Universit\`a di Bologna\\   
piazza di Porta S. Donato 5\\   
I-40126  Bologna, Italy}   
\email{campanin@@dm.unibo.it}   
\thanks{{Partly supported by Italian G.~N.~A.~F.~A, EC grant   
SC1-CT91-0695 and the University of Bologna.   
 Funds for selected   
research topics.}}   
\author{   
Dmitry Ioffe}   
\address{   
Faculty of Industrial Engineering\\   
Technion, Haifa 3200, Israel}   
\email{ieioffe@@ie.technion.ac.il}   
\thanks{Partly supported by the ISRAEL SCIENCE FOUNDATION founded by The Israel   
Academy of Science and Humanities}   
\author{Yvan Velenik}   
\address{   
Laboratoire d'Analyse, Topologie et Probabilit\'es, UMR-CNRS 6632, CMI,   
Universit\'e de Provence, 39 rue Joliot Curie, 13453 Marseille, France}   
\email{velenik\@@cmi.univ-mrs.fr}   
\thanks{{Partly supported by the Swiss National Science Foundation grant  
\#8220-056599.}}   
\vskip 0.2in   
\setcounter{page}{1}   
\begin{abstract}   
We derive precise Ornstein-Zernike asymptotic formula for the decay of the   
 two-point function $\langle \sigma_0\sigma_x\rangle_\beta$   
in the general context of finite range Ising type models on  $\bbZ^d$. The proof   
relies in an essential way on the a-priori knowledge of the  
strict exponential decay of   
the two-point function and,   
by the sharp characterization of phase transition due to Aizenman, Barsky and   
Fern\'{a}ndez,   
goes through in the whole of the high temperature region $\beta<\beta_c$.  
As a byproduct   
we obtain that for every $\beta <\beta_c$, the   
inverse correlation length $\xi_\beta$ is   
an analytic and strictly convex function of direction.    
   
\end{abstract}   
\maketitle   
 
\tableofcontents    
\section{Introduction}   
   
\label{sec_intro}   
The classical Ornstein-Zernike (OZ)  formula \cite{OZ},\cite{Th} gives a sharp  
asymptotic description of  the density pair correlation functions away from the critical  
point.  The original OZ argument is, essentially, a local limit type computation based on  
an ad hoc assumptions on the validity of a certain renewal structure of the correlations.  
   
In this work we give a proof of what happens to be the rigorous counterpart of the OZ   
structural assumption  and, subsequently, derive the corresponding asymptotic formula   
in the context of the high temperature  
finite range Ising models on $\bbZ^d$.  After describing the model   
and formulating the results we give a brief heuristic explanation of the OZ formula in  
terms of the probabilistic local limit theory. In the last part  
 of the Introduction we outline  
the content of the subsequent technical sections.   
   
We would like to mention that although we discuss here only high temperature Ising models,   
many of the ideas we develop could be applied in a broader context of various random  
line type models which possess appropriate uniform   
exponential mixing properties, such  
as, for example, low temperature Pirogov-Sinai interfaces   
in two dimensions.  The   
corresponding study will appear elsewhere.   
  
\subsection{The model}   
\label{ssec_model}   
   
In this work we are considering the class of Ising models with finite-range   
ferromagnetic two-body interactions. To each site $x\in\bbZ^d$ we associate a   
nonnegative real number $J(x) = J(-x) \geq 0$; we suppose that there exists   
$R>0$ such that $J(x)=0$ if $|x|>R$. The collection of these coupling constants   
is denoted by $\uJ$. We consider $\bbZ^d$ as a graph $(\bbZ^d, \cE_\uJ)$, with   
set of vertices $\bbZ^d$ and set of unoriented edges $\cE_\uJ \df   
\setof{(x,y)\in\bbZ^d\times\bbZ^d}{J(x-y)>0}$.   
Let $B \Subset \cE_\uJ$ and $\gb>0$. We denote by $V_B$ the set of sites   
associated to the edges of $B$: given an edge $e\in \cE_\uJ$ and a site   
$x\in\bbZ^d$ we say that $x\in e$ if $x$ is an endpoint of $e$. Then the vertex  
set $V_B$ is defined as,   
$V_B\df \setof{x\in\bbZ^d}{\exists\,e\in B\ {\rm with}\ x\in e}$. The Gibbs measure   
on the graph $(V_B, B)$ at inverse temperature $\gb$ is the probability measure   
on $\{-1,1\}^{V_B}$ defined by   
$$   
\mu_{B,\gb}(\gs) = \frac 1{Z_\gb(B)}\, \exp\bigl\{-\gb \sum_{(x,y)\in B} J(x-y) \gs_x   
\gs_y\bigr\},\quad\quad\quad \gs\in \{-1,1\}^{V_B}\,.   
$$   
A standard argument using Griffiths' second inequality shows that the   
corresponding infinite-volume measure exists; we denote it by $\mu_\gb$.   
Expectation values with respect to the measures $\mu_{B,\gb}$ and $\mu_\gb$, are   
denoted respectively by $\bk{\,\cdot\,}_{B,\gb}$ and $\bk{\,\cdot\,}_\gb$.

\medskip   
   
The central quantity of our study is the 2-point correlation function   
$$   
g_\gb(x)\df\bk{\gs_0\gs_x}_\gb .   
$$   
It plays in the models under consideration precisely the role of the   
density-density correlation function of classical fluids, as can be seen going   
to the lattice gas interpretation of the model, $n_x=\tfrac12 (\gs_x+1)$, where   
the site $x$ is occupied by a particle iff $n_x=1$.

We also introduce the corresponding inverse correlation length: For any   
$x\in\bbR^d$ let   
\begin{equation}  
\label{xi_beta_def}  
\xi_\gb(x) \df - \lim_{k\ra\infty} \frac1{k} \, \log g_\gb([kx]) ,   
\end{equation}  
where for any $y\in\bbR^d$, $[y]\in\bbZ^d$ is the componentwise integer part of   
$y$. A standard sub-additivity argument based on Griffiths' second inequality   
implies that this limit is well-defined and, moreover, letting   
$\vec{\frn}(x)\df x/|x|$,   
\begin{equation}  
\label{griff}  
g_\gb(x) \leq e^{-\xi_\gb(x)} = e^{-\xi_\gb(\vec{\frn}(x))\, |x|}   
\end{equation}  
for all $x\in\bbZ^d$. It also follows from Griffiths' second inequality that the   
function $\xi_\gb$ is convex.

It is important to know for which values of $\gb$ the 2-point function decays   
exponentially, i.e. $\xi_\gb>0$ on $\bbR^d\setminus\{0\}$.  
 Let $\gb_{\rm c}=\gb_{\rm c}(\uJ)$ be the   
inverse critical temperature of the model, i.e.   
$$   
\gb_{\rm c} \df \sup \setof {\gb}{\text{there is a unique Gibbs state at inverse   
temperature $\gb$}}.   
$$   
It is well-known and easy to check that $\infty>\gb_{\rm c}   
 > 0$ when $d\geq 2$. An important   
and highly non-trivial fact is the following theorem due to Aizenman, Barsky and   
Fern\'andez~\cite{ABF}, which asserts that $\xi_\gb$ is in fact an equivalent  
norm on $\bbR^d$,  
\begin{thm}   
$\xi_\gb>0$ if and only if $\gb<\gb_{\rm c}$.   
\end{thm}   
This theorem shows that exponential decay of the 2-point function characterizes   
the high-temperature regime and it provides the basic input for the techniques we develop   
 here.

\subsection{The Results}   
\label{ssec_results}

Our main result describes sharp Ornstein-Zernike-type asymptotics for   
the 2-point function of the models introduced in the previous subsection.   
\begin{bigthm}   
\label{MainTheorem}   
Let $\gb<\gb_{\rm c}$. Uniformly in $\abs{x}\to\infty$   
\begin{equation}   
\label{oz_formula}   
\langle \sigma_0\sigma_x\rangle_\beta~=~\frac{\Phi_\beta\lb   
\vec{\frn}(x)\rb }{\sqrt{\abs{x}^{d-1}}}   
\,{\rm e}^{-\abs{x}\xi_\beta \lb \vec{\frn}(x)\rb}\lb 1\,+\,\so\rb ,   
\end{equation}   
where $\vec{\frn}(x)$ is the unit vector in the direction of $x$;  
$\vec{\frn}(x) = x/\abs{x}$ and $\Phi_\beta$ is a strictly  
positive locally analytic function  
on $\bbS^{d-1}$.   
\end{bigthm}   
   
As a byproduct of the techniques employed for the proof of Theorem~\ref{MainTheorem}  
we deduce that the inverse correlation length $\xi_\beta$ is an analytic and  
strictly convex  
function of the direction.  In order to formulate  
this in a  precise way recall that $\xi_\beta$   
is a convex, homogeneous of order one strictly positive (on $\bbR^d\setminus\lbr 0\rbr$)  
function. As such it is an equivalent norm on $\bbR^d$ and it is the support function of  
the  compact convex set  
\begin{equation}   
\label{intro_Kbeta}   
\Kbeta~=~\bigcap_{n\in\bbS^{d-1}}\setof{t\in\bbR^d}{(t,n)_d\leq \xi_\beta (n)}   
\end{equation}   
with a non-empty interiour $0\in{\rm int}\Kbeta$.  
\begin{bigthm}   
\label{ThmKbeta}   
Let $\beta<\beta_c$. Then $\Kbeta$ has a locally analytic strictly convex boundary  
$\pKbeta$. Furthermore, the Gaussian curvature $\kappa_\beta $ of $\pKbeta$ is   
uniformly positive,   
\begin{equation}   
\label{gaus}   
\bar{\kappa}_\gb\df \min_{t\in\pKbeta}\kappa_\beta (t)~>~0.   
\end{equation}   
\end{bigthm}   
In two dimensions $\Kbeta$ is reminiscent of the Wulff shape (by duality it is   
precisely the low temperature Wulff shape in the  
case of the nearest neighbour interactions).   
The inequality \eqref{gaus} is called then the positive stiffness condition, and one of the  
consequences of Theorem~\ref{ThmKbeta} is the validity of the following strict triangle  
inequality \cite{Io1},\cite{PV2}:  
 Uniformly in $x,y\in\bbR^2$,   
$$   
\xi_\beta (u) +   \xi_\beta (v) - \xi_\beta (u +v)~\geq ~\bar{\kappa}_\gb\lb  
|u| +|v|-|u+v|\rb . 
$$   
In two dimensions $\bar{\kappa}_\gb$ is  
 the minimal radius of curvature of the Wulff shape $\pKbeta$.  
 
In the 2D nearest neighbour case one, using the duality transformation,   
can apply our results   
to  study fluctuations of the $\pm$~interface up to the critical temperature.  
 Similarly, an adjustment of our approach to the general Pirogov-Sinai context in two  
dimensions should, in principle, lead to a   
 comprehensive description  of the fluctuation structure of one-dimensional   
low temperature interfaces.  The corresponding results will appear elsewhere.  
  
Along with $\Kbeta$ we shall consider the set   
\begin{equation}  
\label{Ubeta}  
\Ubeta ~\df ~\setof{x\in\bbR^d}{\xi_\beta (x)\leq 1}~=~  
\setof{x\in\bbR^d}{\max_{t\in\Kbeta}(t,x)_d\leq 1} .  
\end{equation}  
Of course, $\Ubeta$ is just the unit ball in the $\xi_\beta$-norm. It is bounded,  
 convex , and has non-empty interiour for every   
$\beta <\beta_c$.  Furthermore, the polar restatement of Theorem~\ref{ThmKbeta} 
implies that the boundary $\partial \Ubeta$ is also locally analytic and strictly convex
 (c.f. \cite{CIo}). 
 
\medskip   
Results similar or even stronger than    Theorem~\ref{MainTheorem} have been  
obtained in the perturbative  regime $\gb << 1$   
~\cite{AK},\cite{MZ}  and \cite{P-L}. Most recently, the OZ asymptotic has been  
recovered  for the interfaces of the (very) low temperature 2D Blume-Capel  
model in \cite{HK}.   
  
\medskip   
Of course, for the particular case of nearest-neighbour interactions in dimension   
2, the Ornstein-Zernike behaviour of the 2-point function is well-known, through   
explicit computations, see e.g. \cite{MW}. However, non-perturbative,   
dimension independent results of this type have previously been restricted to   
simpler models: Self-avoiding walks have been studied in~\cite{CC}   
(for directions $x$ close to the axis) and~\cite{Io2} (general directions);   
Bernoulli percolation model has been analyzed in~\cite{CCC}   
(for directions $x$ close to the axis) and~\cite{CIo} (general   
directions).   
   
\medskip   
>From a different perspective, Alexander \cite{Al} proved non-perturbative  
lower bounds on two point functions with almost  
the correct order on the prefactor near the decay exponent.   
  Though being  weaker than the sharp asymptotics   
presented in Theorem~\ref{MainTheorem} and failing to capture the fluctuation picture   
behind the phenomenon, these results have the advantage that   
they can be applied to a large variety of models. The core renormalization procedure  
which we   
develop in Section~\ref{Renormalization} is inspired by the ideas of \cite{Al} (see also   
the references therein to his previous works).

\subsection{Probabilistic picture behind the OZ formula}   
\label{ssec_heuristics}   
Let us first explain the order $\norm{x}^{-(d-1)/2}$ of the prefactor in  
\eqref{oz_formula}: Consider a random walk $S_n =V_1+\dots V_n$ on $\bbZ^d$ with  
i.i.d increments $V_i$. Let us assume that the moment generating function  
$\bbE\, {\rm e}^{(t,V_i )_d}$ is finite in a neighbourhood of zero in $\bbR^d$, that the  
distribution of $V_i$ is non-lattice and that the walk $S_n$ is forward in the  
following sense: any point $x$ from the support of the distribution of $V_i$ has  
a positive projection on $\mu\df \bbE V_i$. Given a point $x$ on the direction of  
the principle advance of $S_n$; $\min_{t>0} \norm{x-t\mu}<1$, the probability that   
$S_n$ ``steps'' on $x$ is given by   
\begin{equation}   
\label{sum_iid}   
\sum_{n=1}^{\infty}\bbP\lb S_n =x\rb .   
\end{equation}   
By the usual local limit theorems the term $-\log\bbP\lb S_n =x\rb$ is of the order   
$\norm{x-n\mu}^2/n$. Hence, the main contribution to the above sum comes  
from roughly $\sqrt{\norm{x}}$ terms $n$ around $n_0 = \norm{x}/\norm{\mu}$.   
In other words, up to asymptotically  (with $\norm{x}\to\infty$) negligible terms,   
 the sum in \eqref{sum_iid} is given by the Gaussian summation formula,   
$$   
\frac{c_1}{\sqrt{\norm{x}^d}}\sum_{n=1}^{\infty}{\rm exp}\lbr  
-c_2\frac{(n-n_0)^2}{\norm{x}}\rbr~=~\frac{c_3}{\sqrt{\norm{x}^{d-1}}} .   
$$   
Of course, it is not difficult to give the exact formula  
for $c_3$ in terms of $\mu$ and the covariance   
matrix of $V_i$.   
\vskip 0.1cm  
  
The above sketch almost literally corresponds to the last step  
of the proof of the OZ asymptotic   
formula in the case of the Bernoulli bond percolation in \cite{CIo}. The main effort in the   
latter paper was to show that the percolation cluster from the origin to a (distant) point   
$x\in\bbZ^d$ could be typically split into a density of irreducible pieces with the  
displacements along the endpoints of these pieces playing the role of the i.i.d steps   
$V_1 ,V_2\dots$ of the random walk $S_n$.  
   
In the case of Ising models the two-point function $g_\beta (x)$ also admits a geometric   
random line type representation. Unlike the Bernoulli percolation case, however, different   
portions of this random line interact, whatever splitting rules we employ. In other words   
 in the induced random walk picture the increments $V_1,V_2,\dots$ are dependent. Local   
limit description  of  dependent variables is, in general, a rather delicate matter.  
Fortunately, random lines which show up in the representation of the Ising two-point   
function possess a certain exponential decoupling property. The renormalization procedure   
which we develop in Section~\ref{Renormalization} gives rise to an irreducible splitting   
of the random path in such a way, that the dependence between various sub-paths of the  
splitting has already a uniform exponential decay. The resulting system fits the framework   
of the Ruelle shift operator  on a countable alphabet (of irreducible sub-paths), and,   
as we shall see in the sequel, the associated local limit results are precisely of the  
same analytic nature as in the independent case.

\subsection{Organization of the paper}   
In Section~\ref{Renormalization} we develop  a renormalization procedure leading to  
an irreducible decomposition \eqref{two_point_represent} of the two point function  
$g_\beta (x) = \langle \sigma_0 \sigma_x\rangle_\beta$. This decomposition is translated   
to the Ruelle context in Section~\ref{section_formula}, where we prove both  
Theorem~\ref{MainTheorem} and Theorem~\ref{ThmKbeta}. The proofs rely on general spectral   
properties of the Ruelle operator on countable alphabets described in  
Section~\ref{section_ruelle} and on the local limit analysis of the associated observables   
which is developed in the concluding Section~\ref{section_local}   
\vskip 0.2cm 
\noindent 
{\bf Acknowledgements} D.I. thanks Lev Grinberg for a very careful reading of the  
manuscript and many useful remarks which helped to improve the exposition. M.C. 
and Y.V. gratefully acknowledge the kind hospitality of Technion where part of this work 
was done.

\section{Renormalization}   
\label{Renormalization}   
\setcounter{equation}{0}   
We start by setting up the notation and recalling the well known random   
line representation of the two-point function $g_\beta (x)$. On the  
microscopic level these random paths wiggle in a messy way.  Our main  
renormalization result Theorem~\ref{thm_mass_gap} asserts, however, that on  
sufficiently large scales the random path from $0$ to $x$ exhibits, with an   
overwhelming probability, a regular behaviour, in a sense that it could be   
split into a density of irreducible pieces. The space of irreducible  
paths is defined in Subsection~\ref{sub_irr_paths} and the corresponding   
irreducible representation of $\langle \sigma_0 \sigma_x\rangle_\beta$  
is given by the formula \eqref{two_point_represent} there. The role of   
the cone confinement condition in {\bf (P2)}-{\bf (P4)} will become  
apparent in Section~\ref{section_local}: by the bound \eqref{decoupling}   
it is precisely what one needs in order to represent the system in terms   
of the action of Ruelle operator with a uniformly  
H\"{o}lder continuous potential.   
\subsection{The random-line representation}   
\label{ssection_random_line}   
Recall that given a set of edges   
 $B\Subset\cE_\uJ$ we have defined the associated set of  
vertices as $V_B\df \setof{x\in\bbZ^d}{\exists\,e\in B\ {\rm with}\ x\in e}$.  
For any vertex $x\in V_B$, we define the {\em index} of $x$ in $B$ by  
$\ind x B \df \sum_{e\in B} 1_{e \ni x}$ (as before, $x\in e$ means  
 that $x$ is an  
endpoint of $e$).  
The {\em boundary} of $B$ is defined by $\partial B \df \setof{x\in V_B}{\ind x B  
\text{ is odd}}$.  
  
At each $x\in \bbZ^d$, we fix (in an arbitrary way) an   
ordering  of the $x$-incident edges of the graph:  
$$  
B_x \df \setof{e\in B}{\ind x {\{e\}} >0} =  
\{e^x_1\dots,e^x_{\ind x B}\},  
$$  
 and for two incident edges $e=e_i\in B_x$, $e^\prime = e_j\in B_x$ we say that  
$e\leq e^\prime$ if the corresponding inequality holds for their sub-indices  
; $i\leq j$.   
  
\medskip  
Using the identity $e^{\gb J(e)\gs_t\gs_{t'}} = \cosh(\gb J(e))  
\bigl(1+\gs_t\gs_{t'} \tanh(\gb J(e))\bigr)$, we obtain the following expression  
for the 2-point function of the model in $B$,  
$$  
\bk{\gs_x\gs_y}_{B,\gb} = Z_\gb(B)^{-1} \, \sumtwo{D\subset B}{\partial D =  
\{x,y\}} \prod_{e\in D} \tanh \gb J(e)\,.  
$$  
>From $D\subset B$ with $\partial D = \{x,y\}$, we would like to extract a  
``self-avoiding path''. We use the following procedure:  
\vskip 0.1cm  
\noindent 
\step{1} 
Set $t'_0=y$, $j=0$ and $\Delta_0 =\emptyset$.  
\vskip 0.1cm 
\noindent 
\step{2}  
 Let $e'_j=(t'_j,t'_{j+1})$ be the first edge in  
$B_{t'_j}\setminus\Delta_j$ (in the ordering of   $B_{t'_j}$ fixed above) 
such that $e_j\in D$. This defines $t'_{j+1}$. 
\vskip 0.1cm 
\noindent 
\step{3} Set $\Delta_{j+1}=\Delta_j\cup\setof{e\in B_{t'_j}}{e\leq e'_j}$.  If 
 $t'_{j+1}=x$, then set $n=j+1$ and stop. Otherwise update $j\df j+1$ and return to 
\step{2}. 
\vskip 0.1cm 
This procedure produces a sequence $(t'_0\equiv y,\dots,t'_n\equiv x)$. Let  
$t_k \df t'_{n-k}$ and $e_k \df e'_{n-k}$.  We, thus, constructed a path  
 $\gl\df \gl(D)\df(t_0\equiv x,\dots,t_n\equiv y)$ such that  
\begin{itemize}  
\item $(t_i,t_{i+1})\in B$, $i=0,\dots,n-1$,  
\item $(t_i,t_{i+1}) \neq (t_j,t_{j+1})$ for $i\neq j$,  
\end{itemize}   
(but $t_i=t_j$ for $i\neq j$ is allowed); such a sequence will be called a  
backward edge-self-avoiding line from $x$ to $y$ \footnote{We prefer   
the backward construction of the line $\gl$ because   
it happens to   
be more convenient when reducing to Ruelle's formalism in  
Subsection~\ref{ssec_ReducToRuelle}.}. The construction also yields a set of edges  
$$  
\gD(\gl) \equiv \gD_{n} = \union_{i=1}^{n}  
\setof{e\in B_{t_i}}{e\leq e_i}\,.  
$$  
Notice that $\gD(\gl)$ depends only on $\gl$ (and the order chosen for the edges).  
We use the convenient notation $\sum_{\gl:\,x\mapsto y}$ to represent the  
summation over all (backward) self-avoiding lines from $x$ to $y$. Observe that for any  
$D\subset B$ with $\partial D=\{x,y\}$, $\gl(D)=\gl$ if and only if  
(considering $\gl$ as a set of edges)  
$\gl\subset D$ and $\left(\gD(\gl)\setminus \gl\right) \cap  
D = \emptyset$. We can therefore write  
\begin{equation}  
\bk{\gs_x\gs_y}_{B,\gb} = \sum_{\gl:x\ra y} q_{B,\gb}(\gl)\,,  
\label{eq_RLdef}  
\end{equation}  
where, writing $w(\gl) = \prod_{e\in\gl} \tanh(\beta J(e))$,  
\begin{equation}  
q_{B,\gb}(\gl) = w(\gl)\, \frac {Z_\gb(B\setminus\gD(\gl))}{Z_\gb(B)}\,.  
\label{eq_RLwgt}  
\end{equation}  
Equations \eqref{eq_RLdef} and \eqref{eq_RLwgt} define the {\em random-line  
representation} for the 2-point function of the Ising model on the graph $\cG$.  
It has been studied in detail in \cite{PV1,PV2} and is essentially equivalent  
(though the derivations are quite different) to the  
random-walk representation of~\cite{Aizenman}. We'll need a version of this  
representation on the infinite graph $(\bbZ^d,\cE_\uJ)$. To this end, we use the  
following result (\cite{PV2}, Lemmas~6.3 and 6.9): For all $\beta<\beta_{\rm  
c}$,  
\begin{equation}  
\bk{\gs_x\gs_y}_\beta = \sum_{\gl:x\mapsto y} q_\beta(\gl)\,,  
\label{eq_RLTL}  
\end{equation}  
where $q_\beta(\gl) \df \lim_{B_n \nearrow \cE_\uJ} q_{B_n}(\gl)$ is well defined.  
  
We finally need some rules on how to cut a random-line into pieces.  
Let $\gl=(t_0,t_1,\dots,t_n)$, $z\in\gl$ and let $t_{k(z)}$ be the last hitting  
of $z$ by $\gl$. We write $\gl_<(z) \df (t_0\,\dots,t_{k(z)})$ and $\gl_>(z)  
\df (t_{k(z)},\dots,t_n)$; notice that (as a set of edges) $\gl_<(z) \cap  
\gD(\gl_>(z))=\emptyset$. By the notation $\gl = \gl_1\amalg\gl_2$, we mean  
that there exists $z\in\gl$ such that $\gl_1=\gl_<(z)$ and $\gl_2=\gl_>(z)$.  
Concatenation of 
more than two paths is defined by iterating this procedure, e.g. 
$\lambda_1 \amalg \lambda_2 \amalg \lambda_3 = (\lambda_1 \amalg 
\lambda_2) \amalg \lambda_3$.

We then have the following BK-type inequality:  
\begin{equation}  
\sumtwo{\gl:\,x\ra y}{\gl \ni t} q_\beta(\gl) \leq \sum_{\gl_1:\,x\ra t}  
q_\beta(\gl_1)\,\sum_{\gl_2:\,t\ra y} q_\beta(\gl_2)\,.  
\label{BK}  
\end{equation}  
Indeed, by Griffiths' inequality,  
\begin{align*}   
\sumtwo{\gl:\,x\ra y}{\gl \ni t} q_{\cE,\beta}(\gl) &= \sum_{\gl_2:\,t\ra 
y} q_{\cE,\beta}(\gl_2)   
\sumtwo{\gl:\,x\ra y}{\gl=\gl_1\amalg\gl_2}   
q_{\cE\setminus\gD(\gl_2),\beta}(\gl_1)   
\ =\  \sum_{\gl_2:\,t\ra y} q_{\cE,\beta}(\gl_2)\,   
\bk{\gs_x\gs_t}_{\cE\setminus\gD(\gl_2),\beta}\nonumber\\   
&\leq \sum_{\gl_2:\,t\ra y} q_{\cE,\beta}(\gl_2)\, 
\bk{\gs_x\gs_t}_{\cE,\beta}   
\ =\  \sum_{\gl_1:\,x\ra t} q_{\cE,\beta}(\gl_1)\,\sum_{\gl_2:\,t\ra y} 
q_{\cE,\beta}(\gl_2)\,,   
\end{align*}

\subsection{$K$-skeletons} We coarse-grain microscopic self-avoiding lines  
via an appropriate covering by inflated $\Ubeta$ shapes (see \eqref{Ubeta}):   
Given a self-avoiding line $\lambda = (t_0 ,...,t_n )$  
and a positive number $K>0$ construct the $K$-skeleton $\lambda_K  =(x_0,...,x_N)$  
of $\lambda$   
as follows (Figure~\ref{fig_skeleton}):  
\begin{figure}[t!]  
\input{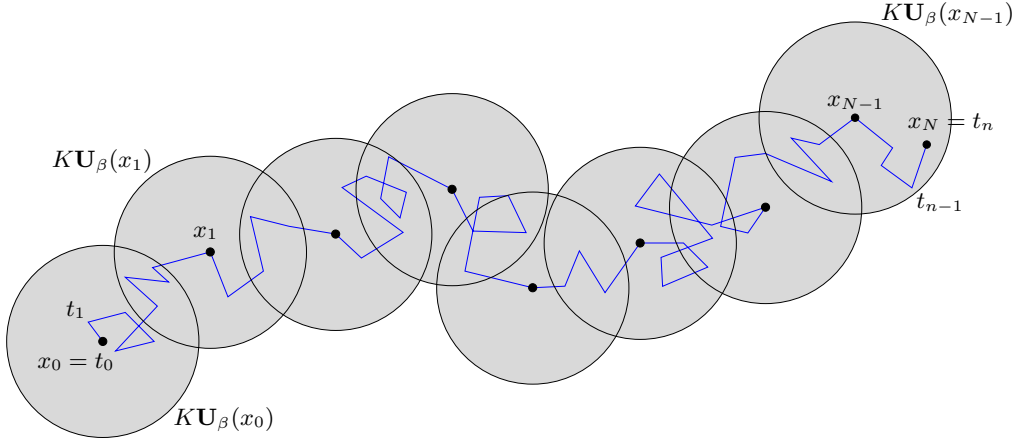}  
\caption{A contour $\gl=(t_0,\dots,t_n)$ and its $K$-skeleton   
$\gl_K=(x_0,\dots,x_N)$.}  
\label{fig_skeleton}  
\end{figure}  
\vskip 0.1cm  
  
\noindent  
\step{1} Set $x_0= t_0$, $j=0$ and $k=0$.  
  
\noindent  
\step{2} If the rest of the line $\{ t_{j+1},...,t_n\}\subseteq K\Ubeta (x_k )$,  
then set $N=k+1$ and $x_N =t_n$ and stop. Otherwise proceed to \step{3}.  
  
\noindent  
\step{3} Find $j^* =\min \setof{i>j}{ t_i\not\in K\Ubeta (x_k )}$.   
Set $x_{k+1}=t_{j^*}$.  
Update $j\df j^*$, $k\df k+1$ and return to \step{2}.  
  
Let us use the notation $\lambda\stackrel{K}{\sim}\lambda_K$ to stress the fact that   
$\lambda_K$ is the $K$-skeleton of $\lambda$. As in the case of paths   
we say that a skeleton $\lambda_K = (x_0 ,...,x_N )$ connects its endpoints,   
$\lambda_K :x_0 \mapsto x_N$.   
  
Of course a particular skeleton  
$\lambda_K =(x_0 ,...,x_N )$ can   
be compatible with many different self-avoiding paths, and we   
introduce the weight  
$$  
\weight{\lambda_K}~=~\sum_{\lambda\stackrel{K}{\sim}\lambda_K}\weight{\lambda} .  
$$  
On any renormalization scale $K$ the BK-inequality \eqref{BK} implies:  
\begin{equation}  
\label{sk_N}  
\weight{\lambda_K}~\leq ~\prod_{l=1}^{N}g_\beta \lb x_l -x_{l-1}\rb ~  
\leq ~{\rm e}^{-(N-1)K}.  
\end{equation}  
  
\subsection{The surcharge inequality}  For $t\in\pKbeta$ let us define  
the surcharge function $\frs_t : {\Bbb Z}^d\mapsto {\Bbb R}_+$ as   
 $\frs_t (x) = \xi_{\beta} (x) - (t,x)_d$.   
Then, given a skeleton $\lambda_K = (x_0 ,...,x_N )$ we   
 define its surcharge as $\frs_t (\lambda_K )=\sum\frs_t (x_{k+1}-x_k )$.   
By the first of   
the inequalities in \eqref{sk_N},  
\begin{equation}  
\label{sk_surcharge}  
\weight{\lambda_K} ~\leq ~{\rm e}^{- (t, x_N )_d - \frs_t (\lambda_K )},  
\end{equation}  
uniformly in $t\in \pKbeta$, scales $K$ and in $K$-skeletons $\lambda_K$.  
 Furthermore, the following crucial surcharge inequality holds:   
\begin{lem}  
\label{lem_surcharge}  
For any $\nu> 0$ there exists a finite renormalization scale   
$K_0 = K_0 (\nu )$ such that    
\begin{equation}  
\label{sk_surcharge_inequality}  
\sumtwo{\lambda_K :0\mapsto x}{\frs_t (\lambda_K )\geq 2\nu |x|}  
\weight{\lambda_K}~\leq ~  
c_1 (\beta ){\rm e}^{-(t, x )_d - \nu |x|} ,  
\end{equation}  
uniformly in $t\in\pKbeta$, $K\geq K_0$ and $x\in {\Bbb Z}^d$  
\end{lem}  
\noindent  
{\em Proof.} There are at most $c_2 (d)K^{d-1}$ choices for each incoming  
skeleton step. Thus, there are at most ${\rm exp}\{ c_3 (d)N\log K\}$ different  
$K$-skeletons of $N$ steps emerging from zero. By \eqref{sk_N} we can   
restrict attention only to those skeletons $\lambda_K :0\mapsto x$ which comprise  
at most $N\leq c_4 |x|/K$ steps.   
Choosing $K_0$ so large that $c_3c_4\log K_0/K_0 <\nu$  
we, in view of the surcharge bound \eqref{sk_surcharge}, arrive to   
the conclusion of the lemma.\qed    
  
There are two types of skeletons with large surcharges which we   
 need to control in order to implement the path decomposition   
procedure:  
  
\subsection{Forward cones and backtracks} Let us fix $\delta \in (0,1/2 )$.   
 For any $t\in\pKbeta$    
define the forward cone  
$$  
\cone{\delta}{t}~=~\setof{x\in {\Bbb Z}^d}{\frs_t (x) <\delta \xi_\beta (x)} .  
$$  
\begin{figure}[t!]  
\input{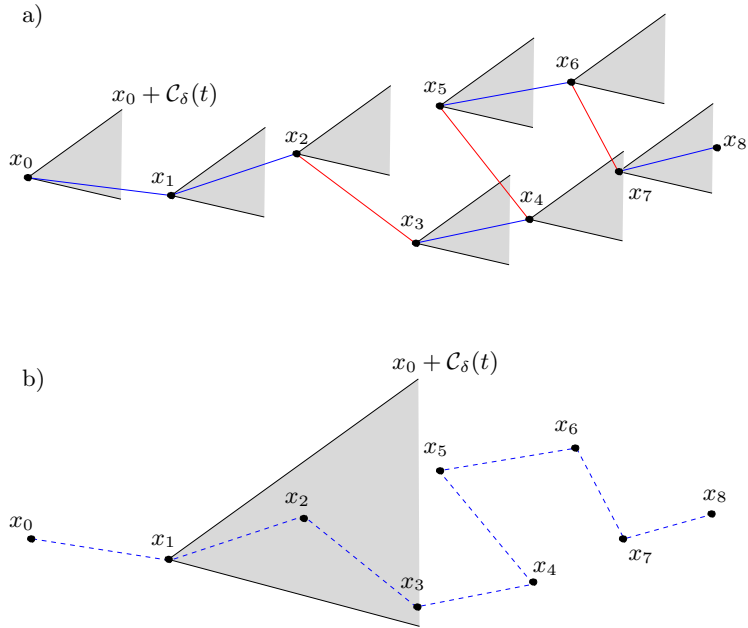}  
\caption{a) A skeleton $\gl_K=(x_0,\dots,x_8)$. The  
increments $[x_2,x_3]$, $[x_4,x_5]$, $[x_6,x_7]$   
are backtracks. Thus, $\nback{t,\delta}(\lambda _K )=3$.  b) The same skeleton   
$\gl_K=(x_0,\dots,x_8)$. The vertex $x_1$ is a cone point of $\gl_K$.}  
\label{fig_points}  
\end{figure}  
Given a $K$-skeleton $\lambda_K = (x_0 ,...,x_N )$ let us define the   
number $\nback{t,\delta}(\lambda _K )$ of $(t,\delta )$-backtracks   
(Figure~\ref{fig_points}~a)) of $\lambda_K$,  
$$  
\nback{t,\delta}(\lambda_K )~=~\#\,\setof{l }{x_{l+1}-x_{l}\not\in\cone{\delta}{t}}.  
$$  
If $x_{l+1} -x_{l}\in \cone{\delta}{t}$,   
  we shall say that  
$x_{l}$ is a forward point of $\lambda_K$.  
  
Notice that the surcharge price of $\lambda_K$ satisfies  
\begin{equation}  
\label{surcharge_back}  
\frs_t \lb \lambda_K\rb~\geq ~\delta K\left(\nback{t,\delta}(\lambda _K ) -1\right) .  
\end{equation}  
  
\subsection{Cone points of skeletons} Given a skeleton $\lambda_K=(x_0 ,...,x_N)$  
let us say that $x_l$ is a $(t,\delta)$-cone point of $\lambda_K$ if   
(Figure~\ref{fig_points}~b))  
$$  
\lbr x_{l+1},...,x_N\rbr~\subset ~x_l +\cone{\delta}{t} .  
$$  
Of course, each cone point of $\lambda_K$ is, in particular, a forward point.  
If a skeleton $\lambda_K$ contains points which do not satisfy the above   
condition, define  
  
\noindent  
$l_1~=~\min\setof{j}{x_j\ \text{is not a}\ (t,\delta)-\text{cone point of}\ \lambda_K}$  
  
\noindent  
$r_1~=~\min\setof{j>l_1}{x_j -x_{l_1}\not\in \cone{\delta}{t}}$  
  
\noindent  
$l_2~=~\min\setof{j\geq r_1}{x_j\   
\text{is not a}\ (t,\delta)-\text{cone point of}\ \lambda_K}$  
  
\noindent  
$r_2~=~\min\setof{j>l_2}{x_j -x_{l_1}\not\in \cone{\delta}{t}}$  
  
\noindent  
$\cdot\, \cdot \,\cdot$  
  
Let us say that $j$ is a $( t, \delta )$-marked point of $\lambda_K$ if it belongs to
the (disjoint) union;    
$j\in\bigvee_k [l_k ,...,r_k )$. Notice that each point of $\lambda_K$ which is not  
marked is, automatically, a $(t,\delta )$-cone (or simply a cone point, if no   
ambiguity with respect to $t$ and $\delta$ arises ) point of $\lambda_K$. We use   
$\nmark{t,\delta}\lb \lambda_K\rb$ to denote the number of all the marked   
points of $\lambda_K$.  
  
\begin{lem}  
\label{lem_mark}  
Uniformly in   
$K$, $\lambda_K$ and $t\in\pKbeta$, the surcharge  
 cost $\frs_{t} (\lambda_K )$ is controlled in terms of the number of   
marked points as   
\begin{equation}  
\label{surcharge_mark}  
\frs_{t} (\lambda_K )~\geq ~\tfrac17\, \delta K\nmark{t,\delta}\lb \lambda_K\rb .  
\end{equation}  
\end{lem}  
  
\noindent  
{\em Proof.} Of course, $\nmark{t,\delta}\lb \lambda_K\rb =\sum_k (r_k -l_k )$.  
We claim that for every marked interval $[l_k,...,r_k )$,  
\begin{equation}  
\label{mark_bound}  
\sum_{j=l_k +1}^{r_k}\frs_t (x_j -x_{j-1})~\geq ~\tfrac17\,\delta K (r_k -l_k ).  
\end{equation}  
Indeed, consider two cases:  
  
\vskip 0.1cm  
\noindent  
\case{1}   
 $(t,x_{r_k} -x_{l_k}) \geq \tfrac27 K (1-\delta )(r_k -l_k)$. Then, since  
 $x_{r_k} -x_{l_k}$ is a back-track, and since $\frs_t$   
evidently inherits  
 from $\xi_\beta$ convexity and homogeneity of order one,   
$$  
\sum_{j=l_k +1}^{r_k}\frs_t (x_j -x_{j-1})  
\geq \frs_t (x_{r_k} -x_{l_k})~  
\geq \delta\, \xi_\beta(x_{r_k} -x_{l_k})  
\geq \delta\, (t,x_{r_k} -x_{l_k})_d  
\geq \delta\, \tfrac17\, K\, (r_k -l_k) .  
$$  
  
\noindent  
\case{2} $(t,x_{r_k} -x_{l_k}) <\frac27 K(1-\delta )(r_k -l_k)$. Notice first  
that $\Kbeta$ is symmetric, so that $t\in\pKbeta \implies -t\in\pKbeta$.  
Therefore, the worst possible displacement of the $t$-projection satisfies  
(recall that we are assuming $K\gg R$)  
$$  
\min_{i\in\{l_k+1,\dots,r_k\}} (t,x_i-x_{i-1}) > -2K\,.  
$$  
This allows us to bound below the number $N_k$ of increments $x_j -x_{j-1}$ from  
the marked interval $j= l_k+1,...,r_k$ that are $(t,\delta)$-backtracking.  
Indeed,  
$$  
(t,x_{r_k} -x_{l_k}) \geq (r_k-l_k-N_k)\, (1-\delta)\, K - N_k\, 2K\,,  
$$  
which gives, since we have fixed the value of $\delta\in (0,1/2 )$,  
$$  
N_k \geq \tfrac17\, (r_k-l_k) .  
$$  
\qed  
  
\subsection{Space of irreducible paths}  
\label{sub_irr_paths}  
 Given $t\in\pKbeta$    
and a path $\lambda= (i_0,\dots ,i_n)$  
let us say that $i_l;~ 0< l <n,$ is a $t$-break point of $\gl$ if $i_l\neq i_k$ for all 
$k\neq l$ and 
$$  
\tilde{\lambda}\cap\lbr i_l +\cH_t\rbr~=~\lbr i_k\rbr ,  
$$  
where $\cH_t$ is the $t$-orthogonal hyper-plane passing through zero; $\cH_t =  
\setof{x\in\bbR^d}{(x,t)_d =0}$, and $\tilde{\lambda}$ is   
the embedding of $\lambda$ with all  
its edges into ${\Bbb R}^d$. Alternatively, $i_l$ is  a $t$-break point of $\gl$ if 
\[ 
\max_{k<l} (i_k ,t)_d \,< \,(i_l,t)_d\,<\,\min_{k<l} (i_k ,t)_d . 
\] 
\begin{figure}[t!]  
\input{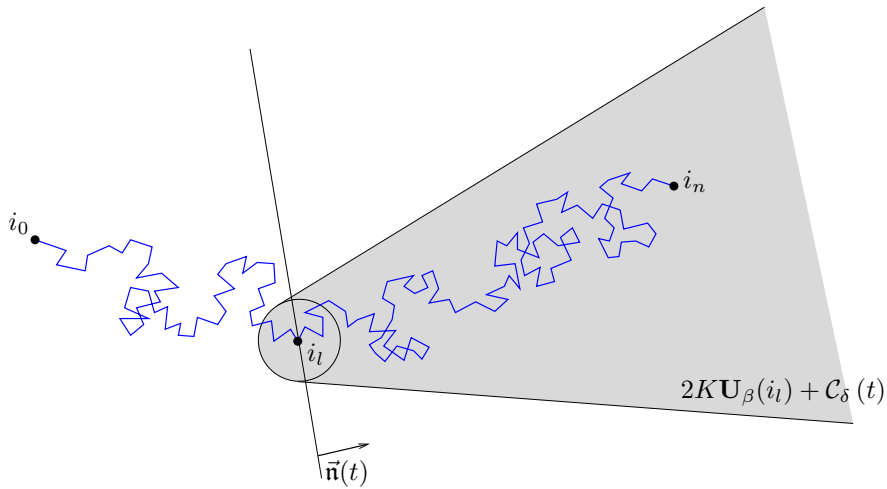}  
\caption{$i_l$ is a $(t,K, \gd )$-correct break point of the contour   
$\gl=(i_0,\dots,i_n)$.}  
\label{fig_c_point}  
\end{figure}  
In addition, given a renormalization skeleton scale $K$ and a forward cone parameter   
$\delta >0$, let us say that a break point $i_l$ of $\gl =(i_0,\dots ,i_n )$ is   
$(t,K, \gd )$-correct if (Figure~\ref{fig_c_point})  
$$  
\lbr i_{l+1},\dots ,i_n\rbr~\subseteq ~ 2K\Ubeta \left(i_l\right)  +\cone{\gd}{t} .   
$$  
  
In particular, if one can find some $(t,\gd )$-cone point  $x_j$ of the  
skeleton $\gl_K$ of $\gl$ such that the break point $i_l$ is on the piece of  
$\gl$ between $x_j$ and $i_n$, and  $i_l\in K\Ubeta (x_j )$, then $i_l$ is  
automatically  $(t,K, \gd )$-correct.   
\begin{thm}  
\label{thm_mass_gap}  
Fix a forward cone parameter $\delta \in (0,1/2 )$. There exist   
a renormalization scale $K_0$ and positive   
numbers $\epsilon =\epsilon (\gd ,\beta )>0$, $\nu =\nu (\gd ,\beta ) >0$ and   
$M=M(\beta)<\infty$,   
such that for all $K\geq K_0$, the upper bound  
 \begin{equation}  
\label{mass_gap}  
\sum_{\gl: 0\mapsto x}\weight{\gl }\1_{\lbr\substack{\text{$\gl$ has less than   
$\epsilon |x|/K$}\\ \text{$(t,K,\delta)$-correct break points}}\rbr}~\leq ~  
M{\rm e}^{-(t,x)_d -\nu |x|} ,  
\end{equation}  
holds   
uniformly in the dual directions $t\in\pKbeta$ and in the end-points $x\in\bbZ^d$.  
\end{thm}  
We relegate the proof of the theorem to the next subsection. Notice, however, that  
 by \eqref{griff} and \eqref{eq_RLdef} the bound \eqref{mass_gap} is trivial 
whenever $t$ and $x$ are such that   
 $x$ lies outside the cone   
$\cone{\nu^{\prime}}{t}$; $\nu^{\prime}=\nu\max_{y\neq 0}|y|/\xi_\beta (y)$.  
  
\begin{figure}[t!]  
\input{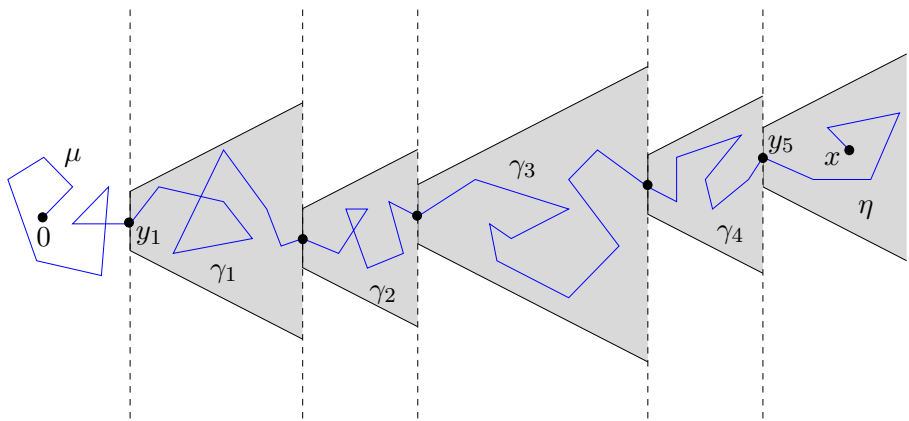}  
\caption{The splitting   
$\gl=\mu\amalg\gga_1\amalg\gga_2\amalg\gga_3\amalg\gga_4\amalg\eta$ into   
irreducible components.}  
\label{fig_split}  
\end{figure}  
For the forward directions $x\in \cone{\nu^{\prime}}{t}$ Theorem~\ref{thm_mass_gap}  
suggests the splitting of a path $\gl :0\mapsto x$ (Figure~\ref{fig_split}):  
\begin{equation}  
\label{path_split}  
\gl ~=~\mu\amalg\gga_1\amalg ,\dots ,\amalg\gga_m\amalg\eta ,  
\end{equation}  
which possesses the following set of properties ${\bf P1}-{\bf P4}$:  
\vskip 0.1cm  
\noindent  
{\bf (P1)}\hskip 0.4cm  All the points  
 $y_1,\dots ,y_{m+1}$ are  break points of $\gl$.  
  
\noindent  
{\bf (P2)}\hskip 0.4cm   
$\eta \subseteq 2K\Ubeta\left(y_{m+1}\right) +\cone{\delta}{t}$ and $\eta$ does  
not contain any $(t,K,\delta)$-correct break point.  
  
\noindent  
{\bf (P3)}\hskip 0.4cm For any $l=1,\dots ,m$, the path $\gga_l$ does not contain   
$(t,K ,\gd)$-correct   
break points, but   
$$  
\gga_l \subset  2K\Ubeta \left( y_{l}\right) +\cone{\delta}{t} .  
$$  
\noindent  
{\bf (P4)}\hskip 0.4cm $\mu$ does not contain   
$(t,K ,\gd)$-correct break points.  
\vskip 0.1cm  
Notice that the successive application of {\bf P1}-{\bf P4} gives an unambiguous   
construction of the decomposition \eqref{path_split}.  
  
Notice, furthermore, that the paths $\gga_l$ (or, more precisely,   
the shifted paths  
 $\gga_l - y_l$) belong to the following basic countable set $S = S(t, K ,\gd) $   
 of  
irreducible paths:  
\vskip 0.2cm  
\noindent  
{\bf Definition}~(The basic set of irreducible paths $S$). Let us say   
that a path $\gga =(i_0 ,...,i_k )\in S$ if  
\begin{enumerate}  
\item \hskip 0.4cm $i_0 =0$ and $(i_0 ,t)< (i_l ,t)< (i_k ,t)$ for all $l=1,\dots ,k-1$.  
\item \hskip 0.4cm $\gga\subset  2K\Ubeta (i_0 ) +\cone{\gd}{t}$.  
\item \hskip 0.4cm $\gga$ does not contain $(t,K, \gd)$-correct break points.  
\end{enumerate}   
Given a path $\lambda =(t_0,\dots ,t_n)$ let us define the displacement   
$V(\lambda)\in\bbZ^d$ along $\lambda$ as the difference between the  
endpoints $V(\lambda )=t_n -t_0$. By Theorem~\ref{thm_mass_gap}, the splitting \eqref{path_split}   
gives rise to the following irreducible representation of the two point   
function $g_\beta (x)= \langle \sigma_0 \sigma_x\rangle_\beta$: Let  
$\epsilon $ be small enough and  
$t\in\pKbeta$ be such that $x\in\cone{\epsilon}{t}$. Then   
\begin{equation}   
\label{two_point_represent}   
g_\beta (x)\left( 1+{\rm o}\left({\rm e}^{-\nu | x |}\right)\right) ~=~   
\sum_{\mu ,\eta}\sum_{m=0}^{\infty}\sumtwo{\gga_1,\dots ,\gga_m \in \cS\,:}{   
V(\mu)+V(\gga_1)+\dots +V(\eta)=x}\weight{\mu\amalg\gga_1\amalg\dots   
\amalg\gga_m\amalg\eta}.   
\end{equation}

 \subsection{Proof of Theorem~\ref{thm_mass_gap}}  
\label{mass_gap_proof}  
The proof is, actually, a modification of the argument developed in \cite{CIo}  
in the context of the Bernoulli bond percolation. It is based on the skeleton  
calculus of the preceding subsections and on the following simple finite energy  
type property:  There exists a positive constant $c_5 >0$, such that for any 
set $B\subset\cE_{\bf J}$ and any path $\gl\subset B$,   
\begin{equation}  
\label{finite_energy}  
q_{B,\beta }(\gl ) 
~\geq ~{\rm e}^{-c_5 \abs{\gga}}.  
\end{equation}  
Notice that all the estimates we employ in the course of the proof   
hold uniformly in $t\in\pKbeta$, and so will the result.    
  
Fix a number $\eta\in (0,1)$.   
\vskip 0.1cm  
\noindent  
{\bf Definition}~    
Given a point $x\in\cone{\delta}{t}$, a skeleton  
scale $K$ and a $K$-skeleton $\gl_K :0\mapsto x$, let us say that $\gl_K$   
is  
$\eta$-admissible if   
 the number of $(t,\delta)$-marked points   
$$  
\nmark{t,\gd}(\gl_K)<\eta\abs{x}/K .  
$$  
By the surcharge inequality \eqref{sk_surcharge_inequality} and the surcharge function   
lower bound \eqref{surcharge_mark}, there exists a   
finite scale $K_0 =K_0 (\eta ,\gd )$, such that   
\begin{equation}  
\label{tune_up1}  
\sumtwo{\gl_K: 0\mapsto x}{\text{$\gl$ is not $\eta$-admissible}}\weight{\gl_K}  
\ \leq~ c_1{\rm exp}\lbr -\frac{\delta\eta}{14}\abs{x} -(t,x)_d\rbr ,  
\end{equation}  
uniformly in the scales $K\geq K_0$ and in $x\in \cone{\gd}{t}$.   
  
Eventually, we are going to pick up $\eta$ sufficiently small, which, as   
the arguments below show, will ensure that up to an  exponentially small   
correction the paths $\gl$ compatible with $\eta$-admissible skeletons contain  
 a density of $(t,K,\gd)$-correct break points, as  has been asserted in   
\eqref{mass_gap} of Theorem~\ref{thm_mass_gap}.  
  
In the sequel we shall   
tacitly assume that the running skeleton scale $K$ is much larger than the   
range of the interaction $R$; $K\gg R$. On every such  skeleton scale $K$   
let us slice $\bbR^d$ into the disjoint  
union of $t$-oriented slabs:   
Let $\vec{\frn}(t)=t/\abs{t}$ be the unit vector in   
the direction of   
$t$  
\begin{equation}  
\label{R_d_slabs}  
\bbR^d~=~\bigvee_{l=-\infty}^{\infty}\lb l\cdot 8K \vec{\frn}(t )~+~\cS_K (t)\rb ,  
\end{equation}  
where  the slab $\cS_K (t)$ is defined via:  
$$  
\cS_K (t)~=~\setof{u\in\bbR^d}{0\leq \left(\frn (t) ,u\right)_d < 8K} .  
$$  
For every $x\in\cone{\gd}{t}$,  
\begin{equation}  
\label{delta_step}  
(t,x)_d~\geq ~(1-\gd )\xi_\beta (x)~\geq ~(1-\gd)\abs{x}\min_{\vec{\frn}\in\bbS^1}  
\xi_\beta (\vec{\frn})~\df ~c_6(\beta )(1-\gd)\abs{x} .  
\end{equation}  
Furthermore, by \eqref{Ubeta} (and in view of the assumption $R\ll K$)   
$(t, x_{k+1}-x_k)_d <2K$,   
whenever $ x_{k+1}-x_k$ is a skeleton increment on the $K$-th skeleton scale.   
 As a result, each skeleton $\gl_K :0\mapsto x$ intersects at least   
$\left[ \frac{c_6 (\beta ) (1-\gd)\abs{x}}{8K}\right]$ subsequent slabs in the  
partition \eqref{R_d_slabs}. On the other hand, if $\gl_K$ is, in addition,   
$\eta$-admissible, then at most $\eta\abs{x}/K$ of these slabs can possibly   
contain marked points of $\gl_K$ .    
The two latter remarks prescribe the choice of the   
number $\eta$:  
\begin{equation}  
\label{eta_choice}   
0~<\eta~<~\frac{(1-\gd)c_6 (\beta )}{16} .  
\end{equation}  
Let us summarize: Given a number $\eta$ as in \eqref{eta_choice} and a skeleton  
parameter $K>K_0 (\eta ,\gd )$, then for any $x\in\cone{\gd}{t}$ and for any   
$\eta$-admissible skeleton $\gl_K : 0\mapsto x$ at least   
$\left[ \frac{c_6 (\beta ) (1-\gd)\abs{x}}{16K}\right]$ of the slabs  
$$  
\cS_{K,l}(t)~\df ~l\cdot 8K \vec{\frn}(t)~+~\cS_K (t)\ ; \qquad l=1,\cdots,   
\left[ \frac{c_6 (\beta ) (1-\gd)\abs{x}}{8K}\right]-1  
$$  
contain only cone points of $\gl_K$. We shall call such slabs $\gl_K$-clean.  
  
 From now on let us fix $\eta$ and $K$ as above. For any $x\in\cone{\gd}{t}$ and  
any $\eta$-admissible skeleton $\gl_K :0\mapsto x$ let us number the   
$\gl_K$-clean slabs in the decomposition \eqref{R_d_slabs} as   
$l_1,\dots ,l_{n}$. As we have just seen,   
\begin{equation}  
\label{n_clean}  
n~=~n (\gl_K )~\geq ~\left[ \frac{c_6 (\beta ) (1-\gd)\abs{x}}{16K}\right],  
\end{equation}  
uniformly in all the situations of interest.   
  
For any $\gl_K$-clean  slab $\cS_{K,l}(t)$ of the skeleton   
$\gl_K =\lb x_0 ,\dots ,x_N\rb$  
let us introduce the indices $i_l$ and $j_l$ via:  
$$  
i_l\,=\,\min\setof{i}{x_i\in\cS_{K,l}(t)}\ \ \ \text{and}\ \ \   
j_l\,=\,\max\setof{j\geq i_l}{x_j\in\cS_{K,l}(t)} .  
$$  
Thus, we can  associate with $\cS_{K,l}(t)$ the embedded sub-skeleton  
$\gl_K^{(l)}=\lb x_{i_l},\dots , x_{j_l}\rb$.   
Similarly,   
let $\gga =\gga_0\amalg\dots\amalg\gga_{N-1}$ be a path compatible with the skeleton $\gl_K$,  
 where $\gga_i: x_{i}\mapsto x_{i+1}$ is the corresponding portion of $\gga$ between the   
skeleton vertices $x_{i}$ and $x_{i+1}$. Then we defined the embedded paths   
$\gga^{(l)}_{-}= \gga_0\amalg\cdots\amalg\gga_{i_l -1}$,    
$\gga^{(l)}=\gga_{i_l}\amalg\dots\amalg\gga_{j_l} $ and   
$\gga^{(l)}_{+}=\gga_{j_l+1}\amalg\dots\amalg\gga_{N-1}$. In this notation,  
$$  
\gga~=~\gga^{(l)}_{-}\amalg\gga^{(l)}\amalg\gga^{(l)}_{+} .  
$$  
\begin{figure}[t!]  
\centerline{\input{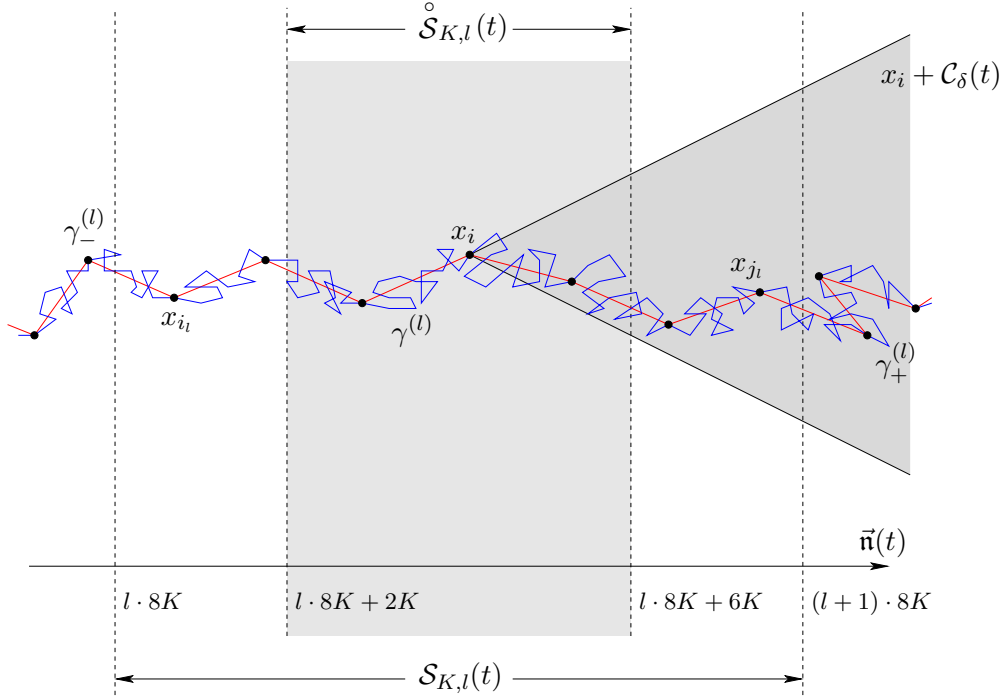}}  
\caption{Clean slab of a skeleton $(x_0,\dots,x_N)$ (only a piece of which is   
drawn): Each point $x_i$; $i=i_l,\dots,j_l$, is a cone point of $\gl_K$. In the   
decomposition $\gga=\gga^{(l)}_-\amalg \gga^{(l)}\amalg\gga^{(l)}_+$ of a path   
$\gga\stackrel{K}{\sim}\gl_K$, the left and right sub-paths $\gga^{(l)}_-$ and   
$\gga^{(l)}_+$ do not intersect $\stackrel{\circ}{\cS}_{K,l}(t)$.}  
\label{fig_clean_slab}  
\end{figure}  
Let us take a closer look at $\gl_K^{ (l)}$ and $ \gga^{(l)}$ (Figure~\ref{fig_clean_slab}):  
 Introducing the inner half-slab  
$$  
\stackrel{\circ}{\cS}_{K,l}(t)~\df~\setof{u\in\bbR^d}{l\cdot 8K +2K\leq (t,u)_d\leq   
 l\cdot 8K +6K},  
$$  
notice that by the very construction   
$x_{i_l},x_{j_l}\in \cS_{K,l}(t)\setminus\stackrel{\circ}{\cS}_{K,l}(t)$.   
In addition, since  
 all the increments of $\gl_K$ on the interval $[i_l,\dots ,j_l)$ are forward;  
$$  
2K~>~\lb x_{i+1}-x_i,t\rb_d~\geq ~(1-\gd )K\qquad \forall~i=i_l,\dots,j_l -1,  
$$  
the  number $j_l -i_l$ of vertices in the sub-skeleton $\gl_K^{(l)}$ is bounded as 
\begin{equation}   
\label{sub_sk_number}  
3\leq j_l -i_l \leq 8/(1-\gd ). 
\end{equation} 
 Finally, the left and right sub-paths  
$\gga^{(l)}_{-}$ and  $\gga^{(l)}_{+}$ are disjoint from   
$\stackrel{\circ}{\cS}_{K,l}(t)$:  
\begin{equation}  
\label{inner_free}  
\gga^{(l)}_{-}\cap \stackrel{\circ}{\cS}_{K,l}(t)~=~\emptyset\qquad\text{and}\qquad  
 \gga^{(l)}_{+}\cap \stackrel{\circ}{\cS}_{K,l}(t)~=~\emptyset .  
\end{equation}  
Consequently, any $t$-break point of   
$\gga^{(l)}$ in the strip $\stackrel{\circ}{\cS}_{K,l}(t)$  
is automatically a $t$-break point of the whole path $\gga$. Furthermore, for any   
$\gl_K$-compatible path $\gga~=~\gga^{(l)}_{-}\amalg\gga^{(l)}\amalg\gga^{(l)}_{+}$  
 one can find $\bar{\gga}^{(l)} :x_{i_l}\mapsto x_{j_l}$, such that   
$\gga^{(l)}_{-}\amalg \bar{\gga}^{(l)}\amalg\gga^{(l)}_{+}$ is still compatible and   
$\gl_K$-compatible, but $\bar{\gga}^{(l)}$ has a $t$-break point in   
$\stackrel{\circ}{\cS}_{K,l}(t)$ and $\abs{\bar{\gga}^{(l)}}\leq c_7 K$. 
By \eqref{sub_sk_number} the total number of all compatible paths is, uniformly in 
$ \gga^{(l)}_{-}, \gga^{(l)}_{+}$ and $\gl_K^{(l)}$, bounded above by $c_8 K^d$.  
Thus, in view of the finite energy condition   
\eqref{finite_energy} applied on the set $B=\cE_{\bf J}\setminus\left( 
\Delta(\gga^{(l)}_{-}) 
 \cup \Delta(\gga^{(l)}_{+})\right)$,  
we infer: 
\begin{equation}  
\label{slab_modification}  
\begin{split}  
\sumtwo{\gga^{(l)}:~\gl_K\stackrel{K}{\sim}\gga^{(l)}_{-}\amalg {\gga}^{(l)}\amalg\gga^{(l)}_{+}}  
{\gga^{(l)}~\text{has no $t$-break points}} \,  
&\weight{\gga^{(l)}_{-}\amalg {\gga}^{(l)}\amalg\gga^{(l)}_{+}}~  
\\  
&\leq ~  
\lb 1-{\rm  e}^{-c_8 K^d}\rb   
\sum_{\gga^{(l)}:~\gl_K\stackrel{K}{\sim}\gga^{(l)}_{-}\amalg {\gga}^{(l)}\amalg\gga^{(l)}_{+}}   
 \weight{\gga^{(l)}_{-}\amalg {\gga}^{(l)}\amalg\gga^{(l)}_{+}}  
\end{split}  
\end{equation}  
The estimate \eqref{slab_modification} is uniform in the points $x\in\cone{\gd}{t}$,   
skeletons $\gl_K :0\mapsto x$, $\gl_K$-clean slabs $\cS_{K,l}(t)$ and in the corresponding   
embedded sub-paths $ \gga^{(l)}_{-}$ and $\gga^{(l)}_{+}$.  Since, by the choice of $\eta$ in   
\eqref{eta_choice} we control the number of different clean slabs of $\eta$-admissible  
 skeletons, \eqref{slab_modification} implies: Let $K\geq K_0 (\eta ,\gd )$. Then   
there exist  $\epsilon =\epsilon (K)>0$  and $\nu =\nu (K)>0$ such that,   
uniformly in $x\in\cone{\gd}{t}$ and in   
the $\eta$-admissible skeletons $\gl_K :0\mapsto x$,   
\begin{equation}  
\label{gl_K-bound}  
\sum_{\gga \stackrel{K}{\sim}\gl_K}\weight{\gga}\1_{\lbr\substack{\text{$\gga$ has less than   
$\epsilon \abs{x}$ break points}\\ \text{inside $\gl_K$-clean slabs}}\rbr}~\leq ~  
{\rm e}^{-\nu\abs{x}}\weight{\gl_K} .  
\end{equation}  
On the other hand any break-point inside a $\gl_K$-clean slab lies inside $K\Ubeta (x_i )$  
 for some cone point $x_i$ of $\gl_K$. Thus, such break points are automatically   
$(t,K,\gd )$-correct, and, thereby,  the claim of Theorem~\ref{thm_mass_gap} follows from   
\eqref{gl_K-bound} and \eqref{tune_up1}. \qed

\section{Ornstein-Zernike formula}   
\label{section_formula}   
The basic decoupling estimate \eqref{decoupling}  which we derive in the  
first subsection, the cone confinement of the irreducible pieces in the  
decomposition \eqref{path_split} and the exponential estimate of Theorem~\ref{thm_mass_gap}   
  enable a reinterpretation of  
the representation formula \eqref{two_point_represent} in terms of the  
Ruelle operator \eqref{L_t_eta} with a uniformly H\"{o}lder continuous \eqref{psi_Holder}  
summable \eqref{mass_gap_ruelle} potential, which paves the way for an application of   
general spectral and  
local limit results of Sections \ref{section_ruelle} and \ref{section_local}.  In view of this   
reinterpretation   
Theorem~\ref{ThmKbeta} more or less directly follows from the analytic perturbation theory   
of non-degenerate eigenvalues as it  
is proved in Subsection~\ref{sub_pKbeta}. The Ornstein-Zernike formula \eqref{oz_formula} is   
derived, along the lines of the general local limit approach of Section~\ref{section_local}, in   
Subsection~\ref{sub_oz_formula}.   
   
\setcounter{equation}{0}   
\subsection{Basic decoupling estimate}   
\label{sub_decouple}  
We prove here an important estimate on the dependence between pieces of a path,  
similar to point 4 of Lemma~5.3 in~\cite{PV1}.  
For two compatible paths $\gga$ and $\gl$ define the conditional  
weight   
$$   
q_{\beta ,\gl}(\gga)~=~\frac{\weight{\gga\amalg\gl}}{\weight{\gl}}.   
$$   
\begin{lem}   
\label{lem_decoupling}   
For every $\beta <\beta_c$ there exists $\theta <1$ and $c_1<\infty$,  
such that for any path $\gga$ and for any pair of compatible and  
 $\gga$-compatible  
paths $\gl = \eta \amalg \gl_1$ and $\gl' = \eta \amalg \gl_2$ with   
$\gD(\gga)\cap( \gD(\gl_1) \cup \gD(\gl_2)) = \emptyset$, the following   
estimate on the ratio of the conditional weights holds:   
\begin{equation}   
\label{decoupling}   
\frac{q_{\beta ,\gl}(\gga)}{q_{\beta ,\gl'}(\gga)}~\geq ~   
{\rm exp}\lbr -c_1\,\sumtwo{t\in\gga}{s\in \gl_1\cup\gl_2}\theta^{|t-s|}\rbr .   
\end{equation}   
\end{lem}   
\begin{rem} 
Lemma~\ref{lem_decoupling} is a principal tool for rewriting the random line weights $q_\gb$ in terms  
of the action of Ruelle operator with H\"{o}lder continuous potential.  In particular, cone 
confinement conditions {\bf (P2)} and {\bf (P3)} have been designed in order to ensure appropriate  
exponential summability properties based on \eqref{decoupling}; see the bound  
\eqref{psi_Holder} below. 
\end{rem}   
\begin{proof}   
Let us consider a finite graph $(\cV,\cE)$ such that $\gga\amalg\gl   
\subset \cE$ and $\gga\amalg\gl' \subset \cE$.   
   
Using $q_{\cE,\beta,\lambda}(\gamma) = 
q_{\cE\setminus\Delta(\lambda),\gb}(\gamma)$, we see that the ratio  
(on the finite   
graph) equals to   
\begin{equation}   
\frac{q_{\cE \setminus\gD(\gl),\gb}(\gga;J)}{q_{\cE\setminus\gD(\gl'),\gb}(\gga;J)}\,.   
\label{ratiorewritten}   
\end{equation}   
However, using the following expression for the weight of a contour (see~(6.39)  
in~\cite{PV2}),   
\begin{equation}   
q_{\cE ,\gb}(\gl;J) = \bigl(\prod_{e\in\gl}\tanh (\gb J(e))\bigr)   
\prodtwo{e\in\gD(\gl)}{e=\bk{t,t'}} \cosh(\gb J(e)) \exp[-\gb J(e) \int_0^1   
\bk{\gs_t\gs_{t'}}^{J_s (e)}_{\cE,\beta}\;\dd s]\,,   
\label{explicitweight}   
\end{equation}   
where $J_s(e) =   
\begin{cases}   
J(e)    & \text{if }e \not\in \gD(\gl)\\   
s J(e)  & \text{if }e \in\gD(\gl)   
\end{cases}   
$, we see that~\eqref{ratiorewritten} is also equal to   
\begin{equation}   
\prod_{e=\bk{t,t'}\in\gD(\gga)} \exp \left[  - \beta J(e) \int_0^1 \left(   
\bk{\gs_t\gs_{t'}}_{\cE\setminus\gD(\gl),\beta}^{J_s} -   
\bk{\gs_t\gs_{t'}}_{\cE\setminus\gD(\gl'),\beta}^{J_s} \right) \dd s   
\right]\,,   
\label{ratioweight}   
\end{equation}   
In view of the strict exponential decay of connectivities in \eqref{griff}, it is   
then sufficient to show that   
$$   
\bk{\gs_t\gs_{t'}}_{\cE\setminus\gD(\gl),\beta}^{J_s} -   
\bk{\gs_t\gs_{t'}}_{\cE\setminus\gD(\gl'),\beta}^{J_s} \leq \sum_{s\in   
\gl_1\cup\gl_2} \bk{\gs_t\gs_{s}}_\beta \bk{\gs_s\gs_{t'}}_\beta \,.   
$$   
Let us prove the latter bound. Let $\cE_1 = \cE\setminus\gD(\gl)$ and $\cE_2 =   
\cE\setminus\gD(\gl ')$.   
We have using the random-line representation:   
\begin{align*}   
\bk{\gs_t\gs_{t'}}^{J_s}_{\cE_1,\gb} &= \sumtwo{\gl:\,t\mapsto   
t'}{\gl\subset \cE_1\cap\cE_2} q_{\cE_1 ,\gb}(\gl;J_s) + \sumtwo{\gl:\,t\mapsto   
t'}{\gl\cap(\cE_1\symdiff\cE_2)\neq\emptyset} q_{\cE_1 ,\gb}(\gl;J_s)\\   
&\leq \sumtwo{\gl:\,t\mapsto t'}{\gl\subset \cE_1\cap\cE_2}   
q_{\cE_1\cap\cE_2 ,\gb}(\gl;J_s) + \sum_{u\in\gl_1\cup\gl_2}   
\sumtwo{\gl:\,t\mapsto t'}{\gl\ni u} q_{\cE_1,\beta}(\gl;J_s)\\   
&\leq \bk{\gs_t\gs_{t'}}^{J_s}_{\cE_1\cap\cE_2,\beta} + \sum_{u\in\gl_1\cup\gl_2}   
\bk{\gs_t\gs_u}^J_{\cE_1,\beta}   
\bk{\gs_{t'}\gs_u}^J_{\cE_1,\beta}\\   
&\leq \bk{\gs_t\gs_{t'}}^{J_s}_{\cE_2,\beta} + \sum_{u\in\gl_1\cup\gl_2}   
\bk{\gs_t\gs_u}_\beta \bk{\gs_{t'}\gs_u}_\beta \,.   
\end{align*}   
The first inequality follows from~\eqref{explicitweight}, 
Griffiths inequality and the fact that all paths containing an edge of  
$\cE_1\symdiff\cE_2$ must also contain a site from $\gl_1\cup\gl_2$; the 
second  
one from the BK-type inequality~\eqref{BK}; finally the last one results 
from another  
application of Griffiths' inequality.  
%
\end{proof}   
\subsection{Reduction to Ruelle's setting}   
\label{ssec_ReducToRuelle}   
Given $t\in\pKbeta$, $\delta \in (0,1)$,  a lattice point  
$x\in\cone{\gd}{t}$ and a path $\gl :0\ra x$ which admits the  
irreducible decomposition \eqref{path_split}, let us rewrite the statistical  
weight $\weight{\gl}$ as   
\begin{equation*}   
\begin{split}   
\weight{\gl}{\rm e}^{(t,x)_d}~   
&=~\weight{\mu \amalg \gga_1 \amalg\dots \amalg \gga_m\amalg   
\eta}{\rm e}^{(t,x)_d}\\   
&=~\weight{\eta}\weight{\mu}{\rm e}^{\lb t,V(\eta )+V(\mu)\rb_d}   
{\rm exp}\lbr \sum_{k=1}^m \psi_\eta^t (\gga_k ,\dots   
,\gga_m )\rbr\,g_{\mu ,\eta}^t (\gga_1,\dots ,\gga_m ),   
\end{split}   
\end{equation*}   
where, as in \eqref{two_point_represent},  we use $V(\gga )$ to denote   
the $\bbZ^d$-displacement between the endpoints of $\gga$ and define  
the potential $\psi_\eta^t$ via:   
\begin{equation}   
\label{psi_eta_t}   
\begin{split}   
{\rm e}^{\psi_\eta^t (\gga_k ,\dots,\gga_m )}~   
&=~\frac{\weight{\gga_k\amalg\gga_{k+1}\amalg \dots\amalg   
\gga_m\amalg\eta}}{   
\weight{\gga_{k+1}\amalg\dots\amalg   
\gga_m\amalg\eta}}\,{\rm e}^{(t,V(\gga_k ))_d}\\   
&= q_{\beta ,\gga_{k+1}\amalg\dots\amalg   
\gga_m\amalg\eta}\left(\gga_k\right) \,{\rm e}^{(t,V(\gga_k ))_d}   
\end{split}   
\end{equation}   
for $k=1,\dots ,m-1$ and, accordingly,  
${\rm e}^{\psi_\eta^t (\gga_m )} =   
 q_{\beta ,\eta}\lb \gga_m\rb  
{\rm e}^{\lb t,V(\gga_m)\rb_d}$.   
   
Similarly, the function $g_{\mu ,\eta}^t$ is defined as   
\begin{equation}   
\label{f_mu_eta_t}   
g_{\mu ,\eta}^t (\gga_1,\dots ,\gga_m ) ~=~   
\frac{\weight{\mu\amalg \gga_1\amalg\dots\amalg   
\gga_m\amalg\eta}}{   
\weight{\gga_{1}\amalg\dots\amalg   
\gga_m\amalg\eta}\weight{\mu}}~=~\frac{q_{\beta ,\gga_{1}\amalg\dots\amalg   
\gga_m\amalg\eta} \lb \mu\rb}{\weight{\mu}}   
\end{equation}   
Notice that since the irreducible paths $\gga_1 ,\dots ,\gga_m$ and the boundary   
condition $\eta$ in the decomposition \eqref{path_split} always satisfy the  
$K$-cone conditions {\bf (P2)} and {\bf (P3)} of Subsection~\ref{sub_irr_paths},  
the decoupling Lemma~\ref{lem_decoupling} implies that the conditional  
weights above are sandwiched between the corresponding unconditional   
ones: There exists $c_2 =c_2 (\theta ,K ,\delta )<\infty $ such that  
\begin{equation}   
\label{free_bounds}   
\frac1{c_2}~\leq ~\frac{q_{\beta ,\gga_{k+1}\amalg\dots\amalg   
\gga_m\amalg\eta}\left(\gga_k\right)}{\weight{\gga_k}} ~\leq ~c_2\qquad{\rm and}\qquad   
\frac1{c_2}~\leq ~ g_{\mu ,\eta}^t (\gga_1,\dots ,\gga_m )~\leq ~c_2  
\end{equation}   
uniformly in $t\in\pKbeta$, $x\in \cone{\gd}{t}$, paths $\gl :0\mapsto x$ and $k=1,\dots ,m$ in the  
decomposition \eqref{path_split} of $\gl$.   
   
In order to enable a uniform local limit study of $g_{\beta }(x)$ along the  
lines of the formalism which will be developed in  
Sections~\ref{section_ruelle} and \ref{section_local} let us, first of  
all, extend any finite sequence of paths $(\gga_1 ,\dots ,\gga_m )$ to an  
infinite one by adding dummy empty paths $\emptyset$. In this way any  
finite sequence of paths $(\gga_1 ,\dots ,\gga_m )$ corresponds to  
the infinite sequence $\ugga =(\gga_1,\dots ,\gga_m ,\emptyset,\emptyset ,   
\dots )$. Thus, given $t\in\pKbeta$ and the forward cone parameter  
 $\delta \in (0,1)$ the basic space $\cS_{\emptyset}$ of infinite sequences   
of irreducible paths can be described as follows:   
\begin{equation}   
\label{space_paths}   
\cS_{\emptyset}~=~   
\bigsetof{\ugga = (\gga_1,\gga_2,\dots )\in\lbr S\cup\emptyset\rbr^{\bbN}}   
{\gga_k=\emptyset\Rightarrow \gga_{j}=\emptyset\ \forall\,j>k} ,   
\end{equation}   
where $S = S(t,\gd )$ is the corresponding space of irreducible paths.   
   
The potential $\psi_\eta^t$ in \eqref{psi_eta_t} has been defined only   
for sequences $\ugga$ of the type  
$\ugga =(\gga_1,\dots ,\gga_m,\emptyset ,\dots )$. However, the basic  
decoupling estimate \eqref{decoupling}   
 implies that for every $t\in\pKbeta$  any two such  
sequences $\ugga$ and $\ugl$ with the proximity  index   
 ${\bf i}(\ugga ,\ugl )\df \min\setof{k}{\gga_k\neq \gl_k}>1$ satisfy the uniform estimate:   
\begin{equation}   
\label{psi_Holder}   
\left| \psi_\eta^t (\ugga )\,- \,\psi_\eta^t (\ugl )\right|~\leq ~c_3   
\theta^{{\bf i}(\ugga ,\ugl )} ,   
\end{equation}   
where the constant $c_3$ depends only on the renormalization scale $K$   
 and on the forward cone parameter $\delta$ which specify the set of   
irreducible paths $S$. Consequently,  $\psi_\eta^t$ admits a  
unique H\"{o}lder continuous extension to the whole of $\cS_{\emptyset}$.   
   
Finally, in view of \eqref{free_bounds},  Theorem~\ref{thm_mass_gap} implies:   
\begin{equation}   
\label{mass_gap_ruelle}   
\sum_{\gga_1\in S}{\rm e}^{\psi_\eta^t (\gga_1 ,\ugga )}~\leq  
~c_4\sum_{x\in\cone{\delta}{t}}{\rm e}^{-\nu |x|}~<~\infty ,   
\end{equation}   
uniformly in $\ugga\in \cS_{\emptyset}$.  
   
   
   
   
   
   
%
   
   
   
   
   
   
   

By \eqref{two_point_represent}  
we have derived the following representation of the two point function:   
 For every   
$x\in \cone{\gd}{t}$;   
   
   
   
%
   
\begin{equation}   
\label{x_representation}   
{\rm e}^{(t,x)_d}g_\beta (x)~=~\smallo\lb {\rm e}^{-\nu |x|}   
\rb  
+\sum_{\mu ,\eta} \weight{\mu}q_\beta (\eta ){\rm e}^{\lb t ,V(\mu )+V(\eta )\rb_d}   
\sum_{n=1}^{\infty}   
\Qn{\mu,\eta}{t}{x-V(\mu )-V(\eta )},   
\end{equation}   
where the weights $\Qn{\mu,\eta}{t}{r}$ are  given by   
\begin{equation}   
\label{Qn_expression}   
 \Qn{\mu,\eta}{t}{r}~=~\sumtwo{\ugga\in\cS_n}{\sum V(\gga_i) =r}   
 {\rm e}^{\Psi^t_{\eta ,n}(\ugga \sep\underline{\emptyset})}g^t_{\mu ,\eta}   
(\ugga ,\underline{\emptyset}),  
\end{equation}   
with  $ \cS_n$ being the set of all $n$-strings of irreducible paths from $S (t,\delta )$, and,  
for every $\ugga\in\cS_n$,   
$$  
{\Psi^t_{\eta ,n}(\ugga \sep\underline{\emptyset})} =   
\psi^t_\eta (\gga_n, \underline{\emptyset})  
+ \psi^t_\eta (\gga_{n-1},\gga_n, \underline{\emptyset})+\dots +  
\psi^t_\eta (\gga_{1},\gga_2 ,\dots ,\gga_n, \underline{\emptyset}) .  
$$  
Thereby,  the weights \eqref{Qn_expression} fall into the general   
framework  of the Ruelle operator induced weights \eqref{Q_n_x}. The local asymptotics   
of the latter are   
studied in general in Section~\ref{section_local}.  
In our case, the associated Ruelle operator $L^t_\eta $ is given by   
\begin{equation}   
\label{L_t_eta}   
L^t_\eta f(\ugga )~=~\sum_{\gga_1\in S}{\rm e}^{\psi^t_\eta (\gga_1 ,\ugga )}   
f(\gga_1 ,\ugga )   
\end{equation}   
By \eqref{psi_Holder} and \eqref{mass_gap_ruelle}  
 $L^t_\eta$ is a bounded linear operator on $\Ggt (\cS_\emptyset )$   
(see Subsection~\ref{sub_setup}   
for the definition of the space $\Ggt (\cS_\emptyset )$ of H\"{o}lder continuous functions on   
$ \cS_\emptyset$  ).   
   
\subsection{The geometry of $\pKbeta$ and the spectral radius $\rho_{\bf S}^t(s )$.}   
\label{sub_pKbeta}   
Since $g_\beta (x)$ is logarithmically asymptotic (see \eqref{xi_beta_def}) to   
${\rm e}^{-\xi_\beta (x)}$,  
$$  
\lim_{|x|\to\infty} \frac1{|x|}\left( \xi_\beta (x) +\log g_\beta (x)\right)\, =\, 0   
$$  
the shape $\Kbeta$ could be alternatively described as the closure of the domain of   
convergence of the series  
$$  
s\,\rightarrow \,\sum_{x\in\mathbb{Z}^d}{\rm e}^{(s,x )_d}g_\beta (x) .  
$$  
Let us fix $t\in \pKbeta$ and $\nu >0$ small. For  
every $|s |<\nu/2$   
 the convergence of the series   
\begin{equation}   
\label{series}   
\sum_{x} g_\beta (x){\rm e}^{(t +s,x)_d }   
\end{equation}   
depends, by the very definition of the surcharge costs,   
  only of the behaviour of $g_\beta (x){\rm e}^{(t +s,x)_d }$ along the   
directions $x$   
satisfying  
$\frs_t (x)\leq \nu |x|$. For such $x$-s, however, the paths $\gl :0\mapsto x$ admit the  
irreducible decomposition \eqref{path_split} with respect to   
the dual direction $t\in\pKbeta$, and we   
are entitled to employ the representation  
\eqref{x_representation}.  Therefore, for $|s| < \nu /2$   
the convergence in  \eqref{series} is  
equivalent to the convergence of the following series:  
\begin{equation}  
\label{series_ruelle}  
\sum_n\sum_{\mu ,\eta } q_\beta (\mu )q_\beta (\eta ){\rm e}^{(t+s ,V(\mu )+ V(\eta ))_d}  
\left[ L_{\eta,s}^t\right]^n g_{\mu ,\eta}^t (\underline{\emptyset}) ,  
\end{equation}  
where we have introduced the ``tilted'' operator  
$$  
L_{\eta ,s}^t f(\ugga )~=~\sum_{\gga_1\in S}{\rm e}^{\psi_\eta^t (\gga_1, \ugga) +\lb  
s ,V(\gga_1 )\rb_d}f(\gga_1 ,   
\ugga )~=~L_\eta^t\lb {\rm e}^{\lb s, V(\cdot )\rb_d}f\rb(\ugga ).   
$$   
\eqref{mass_gap_ruelle} insures that  the operator $L_{\eta ,s}^t$ is well defined for all ${s}<\nu /2$.   
  
By Theorem~\ref{thm_mass_gap} the series  
\begin{equation}  
\label{mu_eta_sum}  
\sum_{\mu ,\eta } q_\beta (\mu )q_\beta (\eta ){\rm e}^{(t+s ,V(\mu )+ V(\eta ))_d}  
\end{equation}  
converges. On the other hand, \eqref{free_bounds} suggests the substitution of  the   
$\left[ L_{\eta,s}^t\right]^n g_{\mu ,\eta}^t (\underline{\emptyset})$   
terms in \eqref{series_ruelle}  
 by   
$$  
\left[ L_{\eta,s}^t\right]^n\1 (\emptyset )\, =\, \sum_{\ugga \in\cS_n}   
{\rm e}^{\Psi^t_{\eta ,n}(\ugga \sep\underline{\emptyset}) +\sum (s, V(\gga_i 
))_d},  
$$  
where $\1$ denotes the constant function on  $\cS_\emptyset$. 
As the in the cases of \eqref{free_bounds} and \eqref{psi_Holder}, the cone confinement   
properties {\bf (P2)} and {\bf (P3)} of the irreducible paths and the basic decoupling  
estimate \eqref{decoupling} imply:   
\begin{equation}  
\label{stringent}  
\sup_n\sup_{\eta ,\eta^\prime}\sup_{\ugga}\,\left\{  
\Psi^t_{\eta ,n}(\ugga \sep\underline{\emptyset}) -  
\Psi^t_{\eta^\prime ,n}(\ugga \sep\underline{\emptyset}) \right\} \,\leq \, c_5 <\infty .  
\end{equation}  
In view of \eqref{mu_eta_sum} this means that the convergence in \eqref{series} is equivalent  
to the convergence of   
$$  
\sum_n \left[ L_{\eta,s}^t\right]^n\1 (\emptyset )   
$$  
for some (and hence for all $\eta$).  
Since \eqref{stringent} evidently implies that   
the spectral radius $\rho_{\bf S}\lb L_{\eta ,s}^t\rb \df \rho_{\bf S}^t (s )$ does not  
depend on $\eta$, we arrive to the following characterization of $\pKbeta$ around $t$: For  
$|s|<\nu /2$,   
\begin{equation}   
\label{eq_pKbeta}   
t+s\in\pKbeta~\Longleftrightarrow~\rho_{\bf S}\lb L_{\eta , s}^t\rb   
\df \rho_{\bf S}^t (s  ) =1 \, .   
\end{equation}   
  
  
Moreover, by Theorem~\ref{thm_mass_gap}, the conditions   
{\bf A1} and {\bf A2} of Section~\ref{section_local} are satisfied for the  
path displacement observable $V :S\mapsto \bbZ^d$.   
Consequently, by the analytic perturbation theory  and  
the non-degeneracy of Hess$(\log\rho_{\bf S}^t)(0)$ established in  
Subsection~\ref{sub_hessian} below, the  equation \eqref{eq_pKbeta} implies that the   
 compact surface $\pKbeta$ is locally analytic and has a uniformly  
positive Gaussian curvature. In particular, the map   
$$   
t~\mapsto ~\frac{\nabla  \rho_{\bf S}^t (0)}{|\nabla \rho_{\bf S}^t (0)|}   
$$   
is a diffeomorphism from $\pKbeta$ to $\bbS^{d-1}$.  
Since by  the general dual description  
of  support functions  $\xi_\gb (x) =(t,x)_d$ if and only if $x$ is orthogonal to a supporting hyperplane 
to $\Kbeta$ at $t$, we conclude: 
For any  
$x\in\bbR^d\setminus 0$ and $t\in\pKbeta$,   
\begin{equation}   
\label{dual_diff}   
\xi_\beta (x) = (t ,x)_d~\Longleftrightarrow ~ \exists\,\alpha\in\bbR_+\  
\text{such that}\  
x=\alpha \nabla\log \rho_{\bf S}^t (0) .   
\end{equation}   
   
\subsection{Proof of the OZ formula}   
\label{sub_oz_formula}   
   

We shall recover the asymptotic behaviour of   
the two point function $g_\beta (x)= \langle \sigma_0\sigma_x\rangle_\beta$ from  
the representation \eqref{x_representation}. The crucial fact is that the  
local limit analysis which will be developed in Section~\ref{section_local} applies for the  
operators $L^t_\eta$ (defined in \eqref{L_t_eta}) and the functions  
$g^t_{\mu ,\eta}$ (defined in \eqref{f_mu_eta_t})  
 {\em uniformly} in $t\in \pKbeta$ and in boundary conditions $\mu$,$\eta$  
satisfying properties {\bf (P2)} and {\bf (P4)} of Subsection~\ref{sub_irr_paths}.   
Indeed, in the language of Section~\ref{section_ruelle} the inequalities \eqref{free_bounds}  
and \eqref{psi_Holder}   
 imply that   
$$   
\sup_{t ,\mu ,\eta}\lbr \snormgt{\psi_\eta^t}+\normgt{g_{\mu ,\eta}^t}+   
\normgt{L^t_\eta}\rbr~<~\infty .   
$$   
In particular (see Theorem~\ref{thm:perron} below), there exists   
$\epsilon >0$ such that the spectrum   
$\Sigma_{\bf S}\lb L^t_\eta\rb$   
 of $L^t_\eta$ satisfies   
$$   
\Sigma_{\bf S}\lb L^t_\eta\rb \bigcap \setof{\mu}{\abs{\mu} >(1-2\epsilon )\rho_{\bf S}^t (0)}~=~   
\lbr \rho_{\bf S}^t (0)\rbr  
$$   
uniformly in $t\in\pKbeta$ and in the boundary conditions $\eta$. Consequently one can  
find an open neighbourhood $\cU$ of the origin in $\bbC^d$, such that the  
family of analytic functions  (see Subsection~\ref{sub_spectral_thm} for the definition of  
the spectral projector  $P_L$), 
$$   
\lbr\,\xi \,\mapsto\,P_{ L^t_{\eta, \xi}}\1 (\underline{\emptyset})\rbr_{t ,\eta}   
$$   
is uniformly continuous on $\overline{\cU}$. By the second of the inequalities   
in \eqref{free_bounds} it follows that the family of the analytic functions   
$$   
\lbr\,\xi\,\mapsto\,\chi_{\mu ,\eta}^t (\xi )\,\df\,P_{ L^t_{\eta,\xi}}g^t_{\mu ,\eta}(   
\underline{\emptyset})\rbr_{t ,\mu , \eta}   
$$   
is uniformly bounded away from zero and infinity on $\cU$. By the Cauchy   
formula the sequence $\lbr \nabla \chi_{\mu ,\eta}^t (0 )\rb$ is also uniformly bounded   
 in $t\in\pKbeta$ and boundary conditions $\mu$ and $\eta$ satisfying  
properties {\bf (P2)} and {\bf (P4)} of Subsection~\ref{sub_irr_paths}.   
\vskip 0.2cm   
  
By the preceding discussion the  asymptotic results of Subsection~\ref{sub_local_setup}  
below   
hold uniformly in $t\in \pKbeta$ and in the boundary  
conditions $\mu, \eta$. For each  
particular choice of the data we shall distinguish between three  
different cases:   
   
Let us fix (see \eqref{R_n_nu} below ) $\nu\in (0,1/2)$ and define   
$$   
R_{n,\nu}^t\,=\,\setof{r\in\bbZ^d}{\abs{r-n\nabla\log\rho_{\bf S}^t (0)} <n^{1-\nu}} .   
$$   
\vskip 0.1cm   
\noindent   
\case{1} $r\in R_{n,\nu}^t$. Then, by Theorem~\ref{thm:local}   
\begin{equation}   
\label{Q_estimate_1}   
\bbQ_{n,\mu ,\eta}^t (r)\, =\, \frac{\chi_{\mu ,\eta}^t (0)}   
{\sqrt{(2\pi n)^d{\rm det}A_{\bf S}^t}}   
{\rm exp}\lbr -\frac1{2n}\cA_{\bf S}^t \lb r-n\nabla\log \rho_{\bf S}^t (0)   
\rb\rbr\lb 1+\so\rb ,   
\end{equation}   
where the quadratic form $ \cA_{\bf S}^t$ is given by  
$\cA_{\bf S}^t (v)= \lb \left[ A_{\bf S}^t\right]^{-1} v,v\rb_d$ and   
$A_{\bf S}^t\df{\rm Hess}\left(\log\rho_{\bf S}^t \right)(0)$.   
\vskip 0.1cm   
   
Pick now a large enough number $M$.   
\vskip 0.1cm   
   
\noindent   
\case{2}  $r\not\in R_{n,\nu}^t$, but $\abs{r}\leq Mn$. Then, as it follows from  
Lemma~\ref{lem:ld_upper},   
\begin{equation}   
\label{Q_estimate_2}   
\bbQ_{n,\mu ,\eta}^t (r)\,\leq {\rm e}^{-c_6 n^{1-2\nu}}\,   
\leq \,{\rm e}^{-c_7 \abs{r}^{1-2\nu}} .   
\end{equation}   
   
\vskip 0.1cm   
\noindent   
\case{3}  Finally, let $r\not\in R_{n,\nu}^t$, and  $\abs{r}> Mn$. In view of  
Theorem~\ref{thm_mass_gap},   
\begin{equation}   
\label{Q_estimate_3}   
\bbQ_{n,\mu ,\eta}^t (r)\,   
\leq \,{\rm e}^{-c_8 \abs{r}} ,   
\end{equation}   
once $M$ has been chosen large enough. This is just an exponential form of   
Markov's inequality.   
\vskip 0.2cm   
Turning back to the expansion \eqref{x_representation}, for each $x\in\bbZ^d$   
define the dual direction $t=t(x)\in\pKbeta$ and the coefficient $\alpha =\alpha (x)$   
as in \eqref{dual_diff}. Set also $n_0 =n_0 (x) =[\alpha (x)]$. Of course,  
\begin{equation}   
\label{x_n_0}   
n_0 (x)~=~\frac{\abs{x}}{\abs{\nabla\log\rho_{\bf S}^t(0)}}\lb 1+\so\rb   
\end{equation}   
uniformly in $\abs{x}\to\infty$.   
   
For every pair of boundary conditions $(\mu ,\eta )$ with  
\begin{equation}   
\label{range_bc}   
\abs{V(\mu )}+\abs{V(\eta )}\,\leq \, n_0^{1/2 -\nu}   
\end{equation}   
we, using the  asymptotic estimates \eqref{Q_estimate_1},   
 \eqref{Q_estimate_2} and  \eqref{Q_estimate_3}, infer that the second sum in   
\eqref{x_representation}  admits the following uniform (in $\abs{x}\to\infty$   
 and in $(\mu ,\eta )$ satisfying \eqref{range_bc})   
asymptotic expression:   
\begin{equation}   
\label{second_sum}   
\frac{\chi_{\mu ,\eta}^t (0)}{\sqrt{(2\pi n_0)^{d-1}   
\cA_{\bf S}^t \lb\nabla\log\rho_{\bf S}^t(0)\rb {\rm det}A_{\bf S}^t}}\lb 1+\so\rb~\df~   
\frac{\phi_{\mu ,\eta}(t)}{\sqrt{\abs{x}^{d-1}}}\lb 1+\so\rb ,   
\end{equation}   
with  
$$   
\phi_{\mu ,\eta}(t)~=~\frac{\chi_{\mu ,\eta}^t (0)\sqrt{\abs{\nabla\log\rho_{\bf S}^t(0)}^{d-1}}   
}{\sqrt{(2\pi )^{d-1}   
\cA_{\bf S}^t \lb\nabla\log\rho_{\bf S}^t(0)\rb {\rm det}A_{\bf S}^t}} .   
$$   
On the other hand, in view of the irreducibility of the boundary conditions   
$(\mu ,\eta)$ in the decomposition \eqref{path_split}, the mass-gap estimate  
\eqref{mass_gap} of Theorem~\ref{thm_mass_gap} implies that  
$$   
\weight{\mu}\weight{\eta}{\rm e}^{\lb t,V(\mu )+V(\eta )\rb_d}~\leq ~   
c_9{\rm e}^{-c_{10}\lb \abs{V(\mu )}+\abs{V(\eta )}\rb } ,   
$$   
uniformly in $t\in \pKbeta$ and in $(\mu ,\eta )$. Consequently, the  
total contribution to the right-hand side of \eqref{x_representation}  
from the terms corresponding to those boundary conditions $(\mu ,\eta )$   
which do not comply with \eqref{range_bc} is at most  
${\rm exp}\lbr -c_{11} \abs{x}^{1/2 -\nu}\rbr$ for some $c_{11} >0$. This  
is negligible as compared to \eqref{second_sum}, and  
 the Ornstein-Zernike   
 formula \eqref{oz_formula} follows with the pre-factor  
$\Phi_\beta \lb\vec{\frn}(x)\rb$ being identified as   
$$   
\Phi_\beta \lb\vec{\frn}(x)\rb~=~\sum_{\mu ,\eta }\phi_{\mu ,\eta }(t(x))  
\, q_\beta (\mu)q_\beta (\eta )  
{\rm e}^{\lb t,V(\mu )+   
V(\eta )\rb_d} .   
$$   
\qed

\section{Ruelle's Perron-Frobenius Theorem for Countable Alphabets}   
\label{section_ruelle}   
The results and the methods of this section are not particularly new.   
 A general treatment of the subshifts on countable  
alphabets could be found in \cite{Br} and in \cite{Sa}. Full shifts are   
studied in the recent preprint \cite{Is} based on the earlier work \cite{CIs}.   
Unfortunately, the setup in the abovementioned papers is different from  
ours and we cannot rely directly on the corresponding techniques therein.   
 In particular, in all these works the authors assumed one or another form  
of irreducibility  
 of the shift, whereas in our context it happens to be natural to permit an   
additional transient class. Thus, for the reader's convenience we prefer   
to formulate the theory in a closed form as we need it here, giving  
exact references whenever possible and providing brief proofs otherwise.   
  
\setcounter{equation}{0}   
\subsection{The Setup}   
\label{sub_setup}   
Let $S$ be  a countable set.   
 We use $\cS_n$ to denote  
the set of $n$-strings $\ux =(x_1,...,x_n)$ of elements of $S$ and  
$\cS$ to denote the set of  
countable $\ux = (x_1,x_2,\dots)$ string of elements of   
$S$. Eventually, we shall study functions defined on set of all finite and infinite  
strings,   
$$   
\cS\bigcup\lb\bigcup_{n=1}^\infty \cS_n\rb .   
$$   
It happens to be convenient to introduce a dummy element $\emptyset$ and define   
\begin{equation}   
\label{S_empty_set}   
\Space\, =\,\setof{\ux\in\lbr S\cup\emptyset\rbr^\bbN}{x_i=\emptyset\Rightarrow   
 x_j=\emptyset\ \forall \,j>i}.   
\end{equation}   
In other words, the infinite strings $\cS\subset \Space$, and for  
every $n\in\bbN$ we extend finite strings from   
$\cS_n$ by attaching to it the infinite sequence   
$\uempty$ of empty elements.  
   
For every $\theta \in (0,1)$ one can define the distance ${\rm d}_\theta$   
on $\Space$ via  
$$   
{\rm d}_\theta (\ux ,\uy )= \theta^{{\bf i}(\ux ,\uy)},   
$$   
where the proximity index between  the strings $\ux\neq \uy$ is given by   
$$   
{\bf i}(\ux ,\uy)~\df ~\min\setof{k}{x_k \neq y_k} .   
$$  
Notice that $\cS$ is a closed subset of $\Space$ in the $ {\rm d}_\theta$ metrics.  
  
Given a function $f: \Space \mapsto \bbC$ and a number  
$k\in\bbN$ define the $k$-th variation of $f$,   
$$   
\var{k}{f}~=~\sup_{\setof{\ux ,\uy}{{\bf i}(\ux ,\uy)\geq k}}\bigabs{f(\ux )   
-f(\uy )} .   
$$   
We say that $f$ is continuous (or more exactly locally uniformly continuous)   
if   
$$   
\lim_{k\to\infty}\var{k}{f}~=~0.   
$$   
The space $\cC =\cC (\Space )$ of bounded continuous functions  
equipped with the usual sup-norm $\normsup{\cdot}$  
 is Banach.   
   
Also, given a number $\theta \in (0,1 )$, we say that $f$ is    
uniformly H\"{o}lder continuous (or, equivalently,  uniformly Lipschitz continuous in the  
${\rm d}_\theta$ metrics of $\Space$) if  
$$   
\snormgt{f}~\df ~\sup_{k>1}\frac{\var{k}{f}}{\theta^k}~<~\infty .   
$$   
Of course, $\snormgt{f}<\infty$ does not imply that $\var{1}{f}<\infty$, and   
hence uniformly H\"{o}lder continuous functions can be  unbounded. However,   
the functional space   
$$   
\Ggt =\Ggt (\Space )\,=\, \setof{f}{f\in\cC\ \text{and}\ \snormgt{f}<\infty}   
$$  
is Banach with respect to the norm  
$\normgt{\cdot} =\normsup{\cdot} +\snormgt{\cdot}$.   
\vskip 0.2cm   
\noindent   
Let a real uniformly H\"{o}lder continuous function $\psi; \,\snormgt{\psi}<\infty$, be  
such that   
\begin{equation}   
\label{L_finite}   
\sup_{\ux\in\Space}\sum_{z\in S}{\rm e}^{\psi (z,\ux)}~ <~\infty   
\end{equation}   
Then the linear operator   
\begin{equation}   
\label{ruelle_operator}   
Lf (\ux)~=~\sum_{z\in S}{\rm e}^{\psi (z,\ux)}f(z,\ux )   
\end{equation}   
is well-defined and bounded on both $\cC$ and $\Ggt$.  
   
Furthermore, given a potential $\psi$ as above and an observable  
$V :S\mapsto\bbZ^d$, the complex operator   
\begin{equation}   
\label{L_itau}   
L_{i\tau}f (\ux )~\df~\sum_{z\in S}{\rm e}^{\psi (z,\ux )  
+i(\tau ,V(z))_d}\,  f(z,\ux ),   
\end{equation}   
where $(\cdot,\cdot)_d$ denotes the scalar product in $\bbC^d$, is also  
defined and bounded on $\Ggt$ and $\cC$ for every $\tau\in [-\pi ,\pi)^d$.  
The original operator $L$ corresponds in the latter notation to $\tau =0$.   
   
\subsection{Spectral properties of $L$ and $L_{i\tau}$} 
\label{sub_spectral_thm}  Given a bounded   
linear operator $T$ on $\Ggt$ let  
  $\Sigma_{{\bf S}}(T)$ and $\Sigma_{{\bf F}}(T)$ denote the  
spectrum and, respectively, the Fredholm spectrum \cite{AKPRS} of $T$.  
We use $\rho_{\bf S}(T)$ and $\rho_{\bf F}(T)$ to denote the corresponding   
spectral radiuses. Any point   
$\gl \in \Sigma_{{\bf S}}(T)\cap\setof{\gl}{\abs{\gl} >\rho_{\bf F}(L)}$   
is an isolated eigenvalue of $T$ (\cite{AKPRS}, Subsection~2.6.12), and  
there exists $\epsilon >0$, such that   
$$   
\setof{\mu }{\abs{\mu -\gl}<2\epsilon}\cap \Sigma_{{\bf S}}~=~\lbr \gl\rbr .   
$$   
Furthermore, for such points $\lambda$ the associated spectral  
projector   
\begin{equation}   
\label{projector_gl}   
P_\gl~=~\frac1{2\pi i}\oint_{\abs{\mu-\gl}=\epsilon}\lb  
\mu I-T\rb^{-1}{\rm d}\mu   
\end{equation}   
is finite dimensional. The dimension of ${\rm Range}\lb P_\gl\rb$ is  
called the algebraic multiplicity of $\gl$. An isolated point  
$\gl_0\in \Sigma_{{\bf S}}$ of algebraic multiplicity  
$1$ called a {\bf non-degenerate eigenvalue } of $L$. There   
is a well-developed analytic perturbation theory of non-degenerate   
eigenvalues, which, in our context,   
 leads to crucial local limit type results. We shall describe   
this in detail in   
Section~\ref{section_local}.

With the above notions in mind let us turn to the spectral  
properties of the operators $L$ and $L_{i\tau}$ which were   
defined  in \eqref{ruelle_operator}   
 and in \eqref{L_itau} respectively.  
   
\begin{thm}   
\label{thm:perron}   
Assume that a uniformly H\"{o}lder continuous real interaction potential $\psi$;  
$\snormgt{\psi}<\infty$, satisfies the summability condition  
\eqref{L_finite}. Then for every $V:S\mapsto\bbZ^d$ and for each  
$\tau\in [-\pi ,\pi )^d$ (in particular for $\tau =0$)   
\begin{equation}   
\label{spectral_gap}  
\rho_{{\bf F}}(L_{i\tau })\,<\,\rho_{\bf S}(L).   
\end{equation}   
Furthermore, $\rho_{{\bf S}}=\rho_{{\bf S}}(L)$ is a non-degenerate  
eigenvalue of $L$  
on $\Ggt$ and the corresponding eigenfunction $h$ is strictly positive;   
\begin{equation}   
\label{h_positive_eigen}   
\inf_{\ux\in\Space}h(\ux )\, >\, 0.   
\end{equation}   
Finally,  the rest of the spectrum of $L$ on $\Ggt$ satisfies  
\begin{equation}   
\label{principle_eigenvalue}   
\sup_{\lambda\in \Sigma_{{\bf S}}   
\setminus \rho_{{\bf S}}}\abs{\gl} ~< ~ \rho_{{\bf S}}.   
\end{equation}   
In particular, there exist $C<\infty$ and  $\gd > 0$,  
such that for any $f\in\Ggt$   
one can find a coefficient $c=c(f)$ satisfying:   
\begin{equation}   
\label{decay_to_h}   
\normgt{ \frac1{\rho_{{\bf S}}^n}L^n f~-~c(f)h}~\leq ~C\normgt{f}(1-\gd)^n .   
\end{equation}   
The above coefficient $c(f)$ satisfies $c(f)h = P_L f$, where we use  
$P_L$ to denote   
the spectral projector \eqref{projector_gl} associated with  
$\gl =\rho_{\bf S}$.   
\end{thm}   
The rest of the section is devoted to the proof of the theorem.   
   
\subsection{Fredholm spectrum}   
\label{sub_gap}   
In this subsection we establish  
the spectral gap assertion \eqref{spectral_gap} of Theorem~\ref{thm:perron}.   
Without loss of the generality we may assume that $\rho_{\bf S}(L)=1$.   
By a version of the Nussbaum's formula \cite{Nus1}, \cite{AKPRS} it  
suffices to show that there exists a compact subset  
$K=K(\tau )$ of $\Ggt$ and a number $n=n(\tau )\in\bbN$, such that   
\begin{equation}   
\label{fredholm_nus}   
\sup_{\normgt{f}\leq 1}\,\inf_{g\in K}\,\normgt{L_{i\tau}^n f -g}\,<\,1.   
\end{equation}   
The $n$-th power of $L$ is given by  
$$   
L^n f(\ux )\,=\,\sum_{\uz\in\cS_n}{\rm e}^{\Psi_n (\uz \sep\ux)}f(\uz ,\ux ),   
$$   
where  
\begin{equation}   
\label{shortcut_psi}   
\Psi_n (\uz \sep\ux)~=~\psi (z_n,\ux)+\psi (z_{n-1},z_n,\ux)+...+   
\psi (z_1,z_2,...z_n,\ux).   
\end{equation}   
It is easy to check  that for every $n\in\bbN$ and for every $\ux,\uy\in\Space$   
\begin{equation}   
\label{Psin_norm}   
\Psi_n (\uz \sep\uy) -\Psi_n (\uz \sep\ux)\leq \gb \theta^{{\bf i}(\ux ,\uy)},   
\end{equation}   
with   
\begin{equation}   
\label{beta_psi}   
\gb =\gb (\psi)~=~\frac{\snormgt{\psi}}{1-\theta} .   
\end{equation}   
\begin{lem}   
\label{apriori_control}   
Assume that $\rho_{\bf S}(L)=1$. Then,   
\begin{equation}   
\label{theta_n_sum}   
\sup_n\,\sum_{\uz\in\cS_n}\normgt{{\rm e}^{\Psi_n\lb\uz\sep \cdot\rb}}~\df   
~M=M(\psi )<\infty   
\end{equation}   
\end{lem}   
\begin{proof}   
By the assumption on $\rho_{\bf S}(L)$,  
 $\inf_{\ux} L^n\1 (\ux )\leq 1$, where $\1$ is the constant function identically equal   
 to $1$. Let us pick $\ux_0$ such that   
$L^n \1 (\ux_0 )\leq 2$. Then,  
for every $\uy\in\Space$, we,using \eqref{Psin_norm}, estimate:   
$$   
{\rm e}^{\Psi_n (\uz\sep\uy )}\,=\, {\rm e}^{\Psi_n (\uz\sep\ux_0)}   
{\rm e}^{\Psi_n (\uz\sep\uy ) -\Psi_n (\uz\sep\ux_0 )}\,\leq  
\,{\rm e}^{\Psi_n (\uz\sep\ux_0)}{\rm e}^{\beta\theta}.   
$$   
Therefore, by the choice of $\ux_0$,   
\begin{equation}   
\label{sum_normsup}   
\sup_n\,\sum_{\uz\in\cS_n}\normsup{{\rm e}^{\Psi_n\lb\uz\sep \cdot\rb}}\,   
\leq \,2{\rm e}^{\beta\theta} .   
\end{equation}   
Moreover, since for every $n$, by \eqref{Psin_norm},   
\begin{equation}   
\label{var_Psi_n}   
\var{k}{{\rm e}^{\Psi_n (\uz\sep \cdot)}}~\leq ~   
\normsup{{\rm e}^{\Psi_n (\uz\sep \cdot~)}} \lb {\rm e}^{\beta\theta^k}-1\rb ,   
\end{equation}   
the sup-norm estimate \eqref{sum_normsup} readily implies the  
conclusion of the Lemma with  
$$   
M(\psi )= 2{\rm e}^{\beta\theta}\lb 1+  
\sup_{t\in (0,1]}\frac{{\rm e}^{\beta t}-1}{t}\rb   
$$   
\end{proof}   
Given $n\in\bbN$, $\tau\in [-\pi ,\pi)^d$ and $\uz\in\cS_n$ set   
$$   
g^n_{\tau,\uz}(\ux ) \,=\,{\rm exp}\lbr \Psi_n (\uz\sep \ux)+i\sum_1^n   
\lb \tau ,V(z_k )\rb_d\rbr .   
$$   
By Lemma~\ref{apriori_control} we (assuming that $\rho_{\bf S}(L)=1$)  
obtain the following estimate:   
\begin{equation}   
\label{sum_gnz}   
\sup_n\,\sum_{\uz\in\cS_n}\normgt{\gnz}~\leq ~M <\infty .   
\end{equation}   
We shall construct the compact sets $K(\tau )$ in \eqref{fredholm_nus} from   
finite linear combinations of functions from the family  
$\{ \gnz\}_{\uz\in\cS_n}$:   
 Given $f\in\Ggt$ with $\normgt{f}\leq 1$  
and $n\in\bbN$ let us represent $L_{i\tau}^{n+1}f$ as   
\begin{equation}   
\label{f_representation}   
L_{i\tau}^{n+1}f (\ux )~=~\sum_{\uz\in\cS_n}\gnz (\ux )L_{\uz ,i\tau}f(\ux ),   
\end{equation}   
where the operator $L_{\uz ,i\tau}$ is defined by   
$$   
L_{\uz ,i\tau}f(\ux )~=~\sum_{u\in S}{\rm e}^{\psi (u,\uz,\ux) +   
i(\tau ,V(u))_d}   
f(u,\uz,\ux ) .   
$$   
Using the obvious inequalities: For every $\phi_1 ,\phi_2\in\Ggt$,  
$\normsup{\phi_1 \phi_2}\leq\normsup{\phi_1}\normsup{\phi_2}$ and, for each  
$k\in\bbN$,   
$$   
\var{k}{\phi_1 \phi_2}~\leq ~\normsup{\phi_1}\var{k}{\phi_2}+   
\normsup{\phi_2}\var{k}{\phi_1},   
$$   
we infer from \eqref{sum_gnz}:   
\begin{equation}   
\label{Litau}   
\normsup{L_{\uz ,i\tau}f(\cdot )}\leq M \ \ \text{and}\ \  
\snormgt{L_{\uz ,i\tau}f(\cdot )}~\leq ~2M\theta^n,   
\end{equation}   
uniformly in $n$, $\uz\in\cS_n$ and in $\normgt{f}\leq 1$. Fix  
now  a large enough power   
 $n$ satisfying $4M^2\theta^{n}<1/2$ and a reference point  
$\ux_0\in\Space$. Defining the coefficients   
$$   
{\bf a }_{\uz}[f]~=~L_{\uz ,i\tau}f (\ux_0 ),   
$$   
we can rewrite \eqref{f_representation} as   
\begin{equation}   
\label{f_coefficients}   
L_{i\tau}^{n+1}f (\ux )~=~\sum_{\uz\in\cS_n}{\bf a }_{\uz}[f]\gnz (\ux )\ +\  
\sum_{\uz\in\cS_n}\gnz (\ux )\lb  L_{\uz ,i\tau}f(\ux )-   
L_{\uz ,i\tau}f(\ux_0 )\rb .   
\end{equation}   
Since we have adjusted the choice of the power $n$ to the estimates in  
\eqref{sum_gnz} and in \eqref{Litau} (notice that the latter also implies   
$\normsup{L_{\uz ,i\tau}f(\cdot )-L_{\uz ,i\tau}f(\ux_0 )}\leq 2M\theta^n$), we   
obtain   
\begin{equation}   
\label{reminder}   
\normgt{\sum_{\uz\in\cS_n}\gnz (\cdot )\lb L_{\uz ,i\tau}f(\cdot )-   
L_{\uz ,i\tau}f(\ux_0 )\rb}\, <\, 1/2.   
\end{equation}   
On the other hand, by the first of the inequalities in \eqref{Litau},  
the sequence of the coefficients $\{{\bf a }_{\uz}[f]\}$ is a  
bounded one; $\abs{{\bf a }_{\uz}[f]}\leq M$. Since by \eqref{sum_gnz}   
for every $\epsilon >0$  
one can choose a finite subset $\cS_{n,\epsilon}\Subset\cS_n$ such that   
$$   
\sum_{\uz\not\in\cS_{n ,\epsilon}}\normgt{\gnz}\,<\,\frac{\epsilon}{M},   
$$   
we are able to derive the following estimate which holds uniformly in  
$\normgt{f}\leq 1$:   
$$   
\normgt{   
L_{i\tau}^{n+1}f (\cdot )  -   
\sum_{\uz\in\cS_{n,\epsilon}}{\bf a }_{\uz}[f]\gnz (\cdot )}   
\, <\,\frac12 +\epsilon .   
$$   
It remains to define the compact $K(\tau )\Subset\Ggt$ as the set of  
all $M$-bounded  
linear combinations of the finite family $\{\gnz\}_{\uz\in \cS_{n,\epsilon}}$;   
$$   
K(\tau )\,\df\,\setof{\sum_{\uz\in \cS_{n,\epsilon}}a_{\uz}\,\gnz (\cdot )}{   
\max_{\uz}\abs{a_{\uz}}\leq M},   
$$   
and the target assertion \eqref{fredholm_nus} follows.

\subsection{The principal eigenfunction of $L$}   
\label{sub_eigenfunction}   
Two main complications we encounter here, as compared to the classical  
setup of subshifts over finite alphabets \cite{Bow},\cite{PP}, are the  
non-compactness of the space $\Space$ and the reducibility of the  
shift $(x_1 ,x_2,\dots )\mapsto (x_2,\dots )$ on $\Space$. The latter   
is merely a nuisance. Nevertheless, it precludes   
an immediate reference to \cite{Sa}, where a non-compact version of  
 Ruelle's Perron-Frobenius theorem has been established in a rather   
general irreducible context.   
   
The results on the existence and strict positivity of the principal  
eigenfunction in the form we need them here, that is  
as asserted in Theorem~\ref{thm:perron}, can be deduced from a  
generalized version of   Krein-Rutman theorem \cite{Nus2} on the  
set-condensing linear maps on cones. However, possibly the most  
transparent way to prove \eqref{h_positive} is to use an approximation   
procedure similar to the one suggested in \cite{CIs}: Let us enumerate   
the elements of $S$ as $x_1 ,x_2,x_3, \dots$ For every $N\in\bbN$ define  
the truncated state space $S^{(N)}=\{ x_1,\dots ,x_N\}$, and, accordingly,   
define the space $\SpaceN$ of countable strings of elements   
from $S^{(N)}\cup\{\emptyset\}$  as in \eqref{S_empty_set}.   
   
For every $\theta\in (0,1)$ $\SpaceN$ is a compact shift-invariant  
subset of $\Space$ in the $d_\theta$-distance (the topology does not   
depend on $\theta$, of course). Let us use $\GgtN$ to denote the restriction   
of $\Ggt$ to $\SpaceN$. Proceeding along these lines, given an interaction   
potential $\psi$ which satisfies the assumptions of Theorem~\ref{thm:perron}   
define the truncated operator $L^{(N)}$ on $\GgtN$,   
$$   
L^{(N)}f(\ux )\,=\,\sum_{z\in S^{(N)}} {\rm e }^{\psi (z,\ux )}f(z,\ux ) .   
$$   
By Lemma~\ref{apriori_control},   
\begin{equation}   
\label{eigen_lim}   
\lim_{N\to\infty} \rho_{\bf S}\lb L^{(N)}\rb\,=\,\rho_{\bf S}\lb L\rb .   
\end{equation}   
On the other hand, despite the reducibility,  
 the arguments of \cite{PP}( pp. 22-24)  directly apply in the  
 $(\GgtN ,L^{(N)} )$-setup above. Consequently, there exists a positive  
eigenfunction $h^{(N)}\in \GgtN$;   
$$   
L^{(N)}h^{(N)}\, =\, \rho_{\bf S}\lb L^{(N)}\rb h^{(N)} ,   
$$   
which, moreover, satisfies the following bound:   
\begin{equation}   
\label{min_max_bound}   
\forall\,\ux,\uy\in\SpaceN\ \ \  h^{(N)}(\ux )\geq  
{\rm e}^{-\beta \theta^{{\bf i}(\ux ,\uy )}}h^{(N)}(\uy ) ,   
\end{equation}   
where the constant $\beta$ has been defined in \eqref{beta_psi}. Notice that   
the estimate \eqref{min_max_bound} holds uniformly in the cutoffs $N$.   
   
It is natural to normalize $h^{(N)} (\uempty )=1$, so that for all $N$;  
${\rm e}^{-\beta}\leq h^{(N)}(\cdot )\leq {\rm e}^{\beta}$. Then for  
every $N\in\bbN$ the restriction to  
$\SpaceN$ of the family $\{ h^{(M)}\}_{M\geq N}$ is bounded in  
$\GgtN$. Using the diagonal procedure, one can extract a subsequence,   
$\lbr h^{(M_k )}\rbr$ which converges in the $\normsup{\cdot}$-norm  
on each of the $\SpaceN$ sets. The limiting function, let us call it $h$,   
is defined on $\cup_N\SpaceN$ and inherits the following properties:   
$$   
h(\uempty )=1\qquad\text{and}\qquad\forall\,\ux,\uy\in\bigcup_N\SpaceN\ \  
h(\ux )\geq {\rm e}^{- \beta \theta^{{\bf i}(\ux ,\uy )} }h(\uy ).   
$$   
Therefore, it can be extended by continuity to the whole of $\Space$,  
and it is straightforward to check from \eqref{eigen_lim} and  
\eqref{theta_n_sum} that the  
extension, which we continue to call $h$, is a strictly positive  
principal eigenvalue of $L$;   
\begin{equation}   
\label{h_positive}   
L h\,=\, \rho_{\bf S}(L)h\qquad\text{and}\qquad  
{\rm e}^{-\beta}\leq h(\cdot )\leq {\rm e}^{\beta} .   
\end{equation}   
This establishes \eqref{h_positive_eigen} of Theorem~\ref{thm:perron}   
  
\subsection{Properties of $\rho_{{\bf S}}(L)$}  
\label{sub_pole}   
In principle it is possible to complete the proof of Theorem~\ref{thm:perron}   
along the lines of   
  \cite{Ru} (Proposition~5.4 on p.90) with necessary adjustments due to   
the fact  
that the unit ball of $\Ggt$ is no longer compact in the space   
of continuous functions $\cC$. Such an approach might also lead to   
 additional complications related to the existence of transient states.   
In fact, the invariant measure in our case will be concentrated on the  
infinite strings of elements from $S$ proper and put zero weight on the  
extended (by $\uempty$) finite strings from $\cS_n$.  Since the latter is   
the main object to be studied in the  application to the  
sharp decay asymptotics of the two-point functions  in the high temperature  
Ising models, 
we  
shall follow a different route:   
   
Using \eqref{h_positive}   
 of the previous subsection   
we can normalize $L$, and,   
apart from the conditions imposed on the interaction $\psi$ in the  
statement of Theorem~\ref{thm:perron},   
 there is no loss of generality to assume that   
\begin{equation}   
\label{L_normalized}   
\rho_{\bf S}(L)=1\ \ \text{and, moreover,}\ \ L\1 (\cdot )\,\equiv\,1 ,   
\end{equation}   
 We need to show:   
\vskip 0.1cm   
\noindent   
{\bf (A)}\ $   
 \lambda=1 \ \text{is the only spectral point of}\  \Sigma_{\bf S}(L)\  
\text{ on the  
spectral circle} \ \ \setof{z\in\bbC}{\abs{z} =1}.$   
\vskip 0.1cm   
\noindent   
{\bf (B)}  
The algebraic multiplicity of   
$\lambda =1$ equals to one, or,  
equivalently,  
 ${\rm Range}\lb P_L\rb$  
  is a one-dimensional  
sub-space spanned by the eigenfunction $ \1 (\cdot )$, where, as in the  
statement of Theorem~\ref{thm:perron},  $P_L$ is the  
spectral projector  \eqref{projector_gl} at the principal  
eigenvalue $\gl =1$.   
   
\vskip 0.1cm   
   
Once {\bf (A)} and {\bf (B)} above are verified, we readily recover   
 the remaining exponential convergence result \eqref{decay_to_h}.  
Indeed, by {\bf (B)}, for every $f\in\Ggt$ there exists a number $c(f)$,  
such that $P_L f =c (f) \1$. On the other hand, the spectral   
radius of $L\lb I-P_L\rb$ is, by {\bf (A)} above, strictly less   
than $1$.   
   
Let $\gl\in\setof{z}{\abs{z}=1}\cap \Sigma_{\bf S} $.  
Since we have already established that $\rho_{\bf F}<1$, it  
follows (\cite{AKPRS}, Subsection~2.6.12)   
 that the eigenspace  
 $\cN (\gl I -L)$ is not empty and finite-dimensional.  
Let $h_\lambda $ be an eigenfunction; $L h_\lambda =\lambda h_\lambda$. By the positivity  
of $L$,  
\begin{equation}  
\label{abs_geq}  
L\left| h_\lambda \right| (\ux ) \geq \left| h_\lambda \right| (\ux )\qquad  
\forall\,\ux\in\Space\, .  
\end{equation}  
Since, for every $\ux\in\Space$ and each $n\in\mathbb{N}$,   
 the probability distribution ${\rm e}^{\Psi_n (\cdot  |\ux )}$ is strictly positive on  
$\cS_n$, we infer from \eqref{abs_geq} that  
\begin{equation}  
\label{sup_ux_n}  
\sup_{\ux\in\Space}\left| h_\lambda \right| (\uz ,\ux ) \,=\,  
\sup_{\ux\in\Space}\left| h_\lambda \right| (\ux )\qquad \forall\,n\ \text{and}\   
\forall\, \uz\in\cS_n \, .  
\end{equation}  
Indeed, taking $\sup_{\ux}$ in both sides of \eqref{abs_geq} certainly does not 
change the ``$\geq$''  sign of the latter inequality. On the other hand,  
 $\sup_{\ux\in\Space}\left| h_\lambda \right| (\uz ,\ux ) \leq  
\sup_{\ux\in\Space}\left| h_\lambda \right| (\ux )$ for any $n$ and $\uz$. 
 
The relation \eqref{sup_ux_n} suggests to consider the restriction of $h_\lambda$ to the  
closed shift invariant    
subset $\cS =\setof{\ux}{x_i\neq  
\emptyset\,\,  \forall\, i}$ of $\Space$:  The H\"{o}lder continuity of   
$h_\lambda$ and the  $n\to\infty$ limit in \eqref{sup_ux_n} readily  
imply:  
$$  
\left| h_\lambda \right| (\cdot )\,\equiv\,   
\sup_{\ux\in\Space}\left| h_\lambda \right| (\ux ) .  
$$  
on $\cS$.  
It is  natural to normalize $h_\lambda$ as $\left| h_\lambda \right|\equiv 1$ on $\cS$.   
But then, given any $\ux\in\cS$,  the function   
$$   
\uz~\rightarrow ~ \frac{h_\gl (\uz ,\ux )}{\gl^n h_\gl (\ux )}   
$$   
is also unimodal on every  $\cS_n$; $n=1,2,\dots$. Since  
$$   
\sum_{\uz\in\cS_n}\text{e}^{\Psi_n (\uz |\ux )}   
\frac{h_\gl (\uz ,\ux )}{\gl^n h_\gl (\ux )} =1,   
$$   
the normalization assumption \eqref{L_normalized}, strict positivity of  
the weights ${\rm e}^{\Psi_n (\uz |\ux )} $ and   
elementary convexity  considerations imply that  
\begin{equation}  
\label{n_trick}  
h_\gl (\uz ,\ux )~=~\gl^n h_\gl (\ux )   
\end{equation}  
for every $\ux \in\cS$, $n\in\bbN$ and $\uz \in\cS_n$. Consequently,  
for every $n\in \bbN$,   
$$   
\sup_{\ux ,\uy\in \cS}\left| h_\lambda (\ux )-h_\lambda (\uy )\right|   
~\leq ~\theta^n \normgt{h_\gl } ,   
$$   
or, in other words, $h_\gl$ is a multiple of $\1$ on $\cS$. In particular,   
\eqref{n_trick} already implies that $\gl =1$,  and  
{\bf (A) } follows. Furthermore, since $L^n h_\lambda (\ux )=h_\lambda (\ux )$ 
for every $\ux\in \cS_\emptyset$ and $h_\lambda$ is H\"{o}lder continuous, we, taking  
the limit $n\to\infty$, readily infer that, actually, $h_\lambda \equiv\1$ on the whole of 
$ \cS_\emptyset$ . 
  
\vskip 0.1cm   
In order to prove {\bf (B)} notice, first of all, that the argument above  
implies that the eigenspace $\cN (I-L )$ is, actually, spanned by $\1$, and,  
hence, $\gl=1$ is a simple eigenvalue.  Now,  
since $\gl=1$ is a Fredholm point; $1 >\rho_{\bf F}(L)$, there exists  
a power $n_0 <\infty$, such that  
$(I-L)^{n_0}f = 0$ for every function $f$ from the range of the projector  
$f\in\text{Range}( P_L)$~(\cite{Ka}, Section III.6.5). If $n_0 =1$, then,  
by the preceding   
remark, we are done. Otherwise, if $n_0 >1$, then for every $f\in\text{Range}( P_L )$   
 there exists a number $d(f)\in\bbC$, such that  
$$   
\lb I-L\rb^{n_0 -1}f~=~d(f)\1 .   
$$   
However, the equation $(I-L)g = d\1$ does not have solutions unless $d=0$.   
Indeed, we may assume that both $d$ and $g$ are real and, in addition, that   
$d$ is non-negative. But, by \eqref{L_normalized},   
$$   
d~= ~\inf_{\ux }(g - Lg)(\ux ) ~\leq ~0.   
$$   
As a result $(I-L)^{n_0 -1 }f = 0$ for every $f\in\text{Range}( P_L )$.  
This reduces the discussion back to the case of $n_0 =1$, and {\bf (B)}   
 follows.   
   
\section{Local Limit Theorem}   
\setcounter{equation}{0}   
\label{section_local}  
We continue to work in the framework and the notation of  
Section~\ref{section_ruelle} and derive strong local limit type results   
associated with the Ruelle's operator $L$. The basic tool is to  
use the spectral theory in order to control the analytic expansions of   
the corresponding log-moment generating functions. We refer to  
 \cite{DS} for a thorough  exposition   
of  the local limit analysis of dependent $\bbZ^d$-valued random  
variables.    See also \cite{AD} where similar results in the CLT region
(and, more generally, in the appropriate  scaling regions for various stable laws) 
have been established for Gibbs-Markov maps. 
\subsection{The setup and the result}   
\label{sub_local_setup}   
Let $V:S\mapsto\bbZ^d$ be an observable and $g\in \Ggt$ is a positive  
function; $\inf_{\ux} g(\ux )>0$. Assuming that the potential   
$\psi ;\ \snormgt{\psi}<\infty ,$ satisfies the  
finiteness assumption \eqref{L_finite}, we associate with  
each $\ux\in\Space$ and every $n\in\bbN$ the weight function $\Qnx$ on $\bbZ^d$ via   
\begin{equation}   
\label{Q_n_x}   
\Qnx (r)~=~\sum_{\uz\in\cS_n:\sum_1^n V(z_i)=r}{\rm e}^{\Psi_n\lb \uz\sep\ux\rb} g(\uz ,\ux ).   
\end{equation}   
Our prime task here is to develop a sharp (as $n\to\infty$) asymptotic formula for  
the weights $\Qnx$. The term ``sharp'' will always mean ``up to zero order terms''.   
An example of such a sharp asymptotic expression is provided by \eqref{decay_to_h}:   
 There exists $c_1>0$, such that  
$$   
\sum_r\Qnx (r)~=~L^n g (\ux )~=~ d_g (\ux )\rho_{\bf S}^n (L)\lb 1~+~   
\smallo ({\rm e}^{-c_1 n}) \rb,   
$$   
where $d_g (\ux)= P_L g (\ux )$.  Since $d_g $ is strictly positive and bounded away  
from zero, there is no loss  
of generality to assume that $\Qnx $ is a probability measure on $\bbZ^d$:   
\begin{equation}   
\label{Q_probability}   
\sum_r\Qnx (r)~=~L^n g (\ux )~\equiv~\rho_{\bf S} (L)~=~1 .   
\end{equation}   
The essential assumptions are, of course, those imposed on the observable $V$:   
\vskip 0.1cm   
\noindent   
{\bf A1.} ${\rm Range}(V)$ generates $\bbZ^d$, in particular $V$ is truly  
$d$-dimensional, in the sense that $\forall\xi\in\bbR^d\setminus 0$,  
the scalar  
product $\lb V(\cdot ),\xi\rb_d\neq const$.

\vskip 0.1cm  
   
\noindent   
{\bf A2.}  There exists $\gd >0$, such that for every $\xi\in\bbR^d$  
with $\norm{\xi}<\gd$   
\begin{equation}   
\label{Cramer_x}   
~\sum_{\uz\in\cS_n}{\rm e}^{\Psi_n (\uz \sep\ux )+\sum_1^n \lb \xi ,V(z_k )\rb_d  }   
~<~\infty .   
\end{equation}   
Notice that by the rigidity bound \eqref{Psin_norm}, the assumption {\bf A2} is   
not sensitive to the choices of boundary condition $\ux$ and powers $n=1,2,...$.  
In particular,   
\begin{equation}   
\label{L_xi}   
L_{\xi}f(\ux )~=~\sum_{z\in S}{\rm e}^{\psi (z , \ux ) + (\xi , V(z))_d} f(z,\ux ),   
\end{equation}   
is a well defined bounded linear operator on $\Ggt$ for every $\xi\in\bbC^d $   
 with $\abs{{\rm Re}(\xi )} <\gd$.  
\vskip 0.1cm   
   
Our first result is a rough Gaussian  
large deviation upper bound which enables to focus the  
attention on the values of $r$ near the running average   
$$   
nv_{n,\ux}~\df~\Enx\sum_{k=1}^n V(z_k ) .   
$$   
\begin{lem}   
\label{lem:ld_upper}   
For every $\nu >0$ there exist $c_2,c_3 >0$, such that   
\begin{equation}   
\label{nu_ld}   
\sum_{r:\abs{r-nv_{n,\ux}}\geq n^{1-\nu}}\Qnx (r)~<c_2{\rm e}^{-c_3n^{1-2\nu}}.   
\end{equation}   
\end{lem}   
Lemma \ref{lem:ld_upper} is a standard consequence of the exponential  
Markov inequality and the non-degeneracy condition \eqref{H_n_hessian} which   
is formulated below (and, subsequently, is proved in  
Subsection~\ref{sub_hessian}).

>From now on we fix $\nu\in (0, 1/2 )$ and concentrate on deriving uniform  
sharp asymptotics   
of $\Qnx (r)$ over the set   
\begin{equation}   
\label{R_n_nu}   
R_{n,\nu}~=~\bigsetof{r\in\bbZ^d}{\abs{r-nv_{n,\ux}}< n^{1-\nu}} .   
\end{equation}   
It is exactly on this stage that we shall extensively rely on the spectral   
analysis of Section~\ref{section_ruelle}. In order to structure our main result   
here in an optimal way let us formulate it in the form of several separate   
 propositions:   
\vskip 0.1cm   
   
We claim that there exists an open neighbourhood $\cU$ of the origin in $\bbC^d$,  
such that, uniformly in boundary conditions $\ux\in\Space$, all  the  
properties listed below hold:

\begin{lem}   
\label{lem:perturbation}   
The functions  
$$   
\rho_{\bf S}(\xi)~\df ~\rho_{\bf S}(L_\xi)\qquad{\rm and}\qquad  
\chi_{\ux}(\xi )~\df ~P_{L_\xi } g (\ux )   
$$   
are analytic and bounded away from zero on $\cU$. Furthermore, for every  
$\xi \in \cU$,  $\rho_{\bf S}(\xi )$ is (c.f. Subsection~\ref{sub_setup})   
a non-degenerate eigenvalue of $L_\xi$ and, independently of a particular  
choice of $\xi\in\cU$, there exists $\epsilon >0$ such that the rest of the spectrum  
of $L_\xi$ lies inside the circle of the radius  
$(1-\epsilon)\abs{\rho_{\bf S}(\xi )}$.  
\end{lem}   
This is a rather standard assertion of the analytic perturbation theory based  
on  Theorem~\ref{thm:perron} and assumption {\bf A2}. We shall explain it   
in more detail (and with the appropriate references to \cite{Ka}) in  
Subsection~\ref{sub_perturbation}.   
\vskip 0.1cm   
As it follows from Lemma~\ref{lem:perturbation}, the log-Laplace transforms   
\begin{equation}   
\label{H_n_def}   
\Hnx (\xi)~\df ~\frac1n\log L_\xi^n g (\ux )   
\end{equation}   
are defined and analytic on $\cU$. Moreover,   
\begin{lem}   
\label{lem:H_n_expansion}   
There exist $c_4 >0$ such that, uniformly in $\ux\in\Space$ and $\xi\in\cU$,   
\begin{equation}   
\label{H_n_expansion}   
\Hnx (\xi)~=~\log\rho_{\bf S}(\xi) +\frac1n\log \chi_{\ux}(\xi)  
+\smallo \lb{\rm e}^{-c_4 n}\rb .   
\end{equation}   
In addition the Hessians ${\rm Hess}\lb \Hnx \rb$ are uniformly non-degenerate   
at $\xi =0$;   
\begin{equation}   
\label{H_n_hessian}   
\inf_{\ux\in\Space}\abs{{\rm det}\lb{\rm Hess}(\Hnx )(0)\rb} ~>~0.   
\end{equation}   
\end{lem}   
The non-degeneracy condition \eqref{H_n_hessian} is responsible for the  
Gaussian form of our main uniform local limit result: Define   
$$   
A_{\bf S}  ~=~{\rm Hess}\lb \log \rho_{\bf S}\rb(0) .   
$$   
\begin{thm}   
\label{thm:local} Uniformly in $r\in R_{n,\nu}$ (see \eqref{R_n_nu}) and  
$\ux\in\Space$   
\begin{equation}   
\label{Q_n_asymptotics}   
\Qnx (r)~=~\frac{d_g (\ux )\rho_{\bf S}^n(0)}{\sqrt{(2\pi n)^d{\rm det}(A_{\bf S})}}   
{\rm exp}\left\{-\frac1{2n}\lb A_{\bf S}^{-1}(r-nv_{n,\ux} ),   
(r-nv_{n,\ux} )\rb_d\right\} \lb 1~+~\so\rb .   
\end{equation}   
\end{thm}   
 Notice that since the running average $v_{n,\ux}=\nabla \Hnx (0 )$,  
the uniform analytic  
expansion \eqref{H_n_expansion} implies,   
   
\begin{equation}   
\label{v_n_expansion}   
v_{n,\ux}~=~\nabla\log\rho_{\bf S} (0)+\frac1n\nabla\log\chi_{\ux} (0) +   
\smallo\lb{\rm e}^{-c_5 n}\rb ,   
\end{equation}   
and, consequently, we  
 could have written the $\ux$-independent term  
$\nabla \log\rho_{\bf S}(0)$ instead of $v_{n,\ux }$ in the target  
asymptotic formula \eqref{Q_n_asymptotics}.   
\vskip 0.1cm   
\noindent   
{\em Proof of Theorem~\ref{thm:local}.} The proof is a blend of the conventional   
 local CLT techniques and the (equally conventional) change of measure  
by exponential tilts argument reinforced with an analytic control over log-Laplace   
transforms through the expansion \eqref{H_n_expansion}. We shall merely sketch it  
here with an emphasis on how the spectral analysis of the Ruelle's operator  
enters the picture. We refer to \cite{DS} for a comprehensive general  
exposition of the local limit  
theory and also to \cite{PP}, where similar results are obtained for the  
Ruelle's operators over finite alphabets.    
   
As before, there is no loss of generality to assume that $\Qnx$ is a   
 probability measure on $\bbZ^d$, in particular, that $d_g (\ux )\equiv 1$ and  
that $\rho_{\bf S}(L) =\rho_{\bf S}(0) =1$.  
\vskip 0.1cm   
\noindent   
\step{1} Fix a small $\epsilon >0$. We shall start by proving  
\eqref{Q_n_asymptotics} for the values of $r$ satisfying  
(see \eqref{v_n_expansion})   
\begin{equation}   
\label{close_r_values}   
\abs{r-n\nabla\log \rho_{\bf S}(0)}~\leq ~n^{1/2 -2\epsilon} .   
\end{equation}   
In this case the target asymptotic expression \eqref{Q_n_asymptotics}   
 of Theorem~\ref{thm:local}   
takes a simpler form:   
\begin{equation}   
\label{Q_n_r_close}   
\Qnx (r)~=~\frac1{\sqrt{(2\pi n)^d{\rm det}(A_{\bf S})}}   
 \lb 1~+~\so\rb .   
\end{equation}   
Let $\widehat{\bbQ}_{n,\ux}$ denote the Fourier transform of $\Qnx$,   
$$   
\widehat{\bbQ}_{n,\ux}(\tau )~=~\sum_{t\in\bbZ^d}\Qnx (t){\rm e}^{i (\tau ,t)_d}~=~   
L_{i\tau}^n g (\ux ).   
$$   
By the Fourier inversion formula,   
\begin{equation}   
\label{fourier_formula}   
\Qnx (r)~=~\frac{1}{(2\pi )^d}\,\int ...   
\int_{[-\pi ,\pi]^d}{\rm e}^{-i(\tau ,r)_d}   
\chnx (\tau ){\rm d}\tau .   
\end{equation}   
Given $\gd >0$,  
we split $[-\pi ,\pi ]^d$ into three disjoint regions of integration:   
\begin{equation}   
\label{regions}   
\begin{split}   
[-\pi ,\pi ]^d ~&=~A_\epsilon\vee A_{\epsilon ,\gd}\vee A_\gd~\\   
&\qquad\df ~   
\setof{\tau}{\abs{\tau}<n^{-1/2+\epsilon}}\vee\setof{\tau}{n^{-1/2+\epsilon}   
\leq \abs{\tau}   
<\delta}\vee\setof{\tau}{\abs{\tau}\geq\gd} .   
\end{split}   
\end{equation}   
The integral over $A_\gd$ could be ignored by the virtue of the following  
proposition, which we shall prove in Subsection~\ref{sub_decay_away}:   
\begin{pro}  
\label{pro:off_gap} For every $\gd >0$ there exists $\eta =\eta (\gd )>0$, such that   
\begin{equation}   
\label{off_gap}   
\sup_{\tau\in [-\pi ,\pi]^d\setminus (-\gd,\gd)^d}   
\rho_{\bf S} \lb L_{i\tau}\rb~\leq ~1-\eta .   
\end{equation}   
\end{pro}   
An immediate consequence is that, uniformly in $\ux\in\Space$ and $\tau\in A_\gd$,  
\begin{equation}   
\label{A_delta}   
\abs{\chnx (\tau)}\leq {\rm e}^{n\log (1-\eta )} .   
\end{equation}   
   
\vskip 0.1cm   
\noindent   
Turning to $A_{\epsilon,\gd}$ notice that  
if $\gd> 0$ is sufficiently small, then $iA_{\epsilon ,\gd}\subset\cU$, and,  
consequently,   
$$   
\chnx (\tau)~=~{\rm e}^{n\Hnx (i\tau )} .   
$$   
Choosing, if necessary, $\gd >0$ even smaller, we infer from the analytic  
expansion formula \eqref{H_n_expansion} and the non-degeneracy of   
 ${\rm Hess}\lb\log \rho_{\bf S}\rb (0)$, that there exists $c_5 >0$, such that ,   
\begin{equation}   
\label{A_epsilon_delta}   
\abs{\chnx (\tau )}~\leq ~{\rm e}^{-c_5 n \abs{\tau}^2}~\leq ~   
{\rm e}^{-c_5 n^{2\epsilon}},   
\end{equation}   
uniformly in $\ux\in\Space$ and $\tau \in A_{\epsilon ,\gd}$.   
\vskip 0.1cm   
\noindent   
Finally, uniformly in $\tau$ from the remaining region $A_\epsilon$,   
\begin{equation*}   
\begin{split}   
\chnx (\tau ){\rm e}^{-i(\tau ,r)_d}~&\stackrel{\eqref{close_r_values}}{=}~   
{\rm exp}\lbr n\Hnx (i\tau ) -i n \lb\tau ,\nabla \log \rho_{\bf S}(0)\rb_d +\so\rbr\\   
&\stackrel{\eqref{H_n_expansion}}{=}~   
{\rm exp}\lbr -\frac{n}{2} \lb A_{\bf S}\tau ,\tau \rb_d +\so\rbr,   
\end{split}   
\end{equation*}   
and \eqref{Q_n_r_close} follows.   
\vskip 0.1cm  
\noindent   
\step{2} In order to extend the result to the full range of $r\in R_{n,\nu}$   
as it has been asserted in Theorem~\ref{thm:local}, consider the family of  
``tilted'' measures $\{\Qnxi{\xi}\}$ (indexed by $\xi\in\cU\cap\bbR^d$):   
$$   
\Qnxi{\xi} (r)~=~ \frac{ {\rm e}^{(\xi , r)_d}}{L_\xi g (\ux )}\Qnx (r)~=~   
{\rm exp}\lbr (\xi ,r)_d -n\Hnx (\xi )\rbr\Qnx (r) .   
$$   
The expectation $nv_{n,\ux}(\xi )$ under the measure $\Qnxi {\xi }$ is, according to   
 \eqref{H_n_expansion}, given by the following asymptotic expression:   
$$   
v_{n,\ux}(\xi )~=~\nabla \log \rho_{\bf S}(\xi )+\frac1n\nabla \log\chi_{\ux }(\xi )   
+\smallo\lb{\rm e}^{-c_6 n}\rb .   
$$   
Since the Hessian ${\rm Hess}(\log \rho_{\bf S})$ is non-degenerate at $\xi=0$, we,   
actually independently from $\ux\in\Space$,    
can pick a small  $\gd >0$, such that the map $\xi\mapsto v_{n,\ux }(\xi )$ has  
an analytic inverse on $\setof{\xi}{\abs{\xi} <\gd}\subset \cU$. Since, in this   
case,   
$$   
R_{n,\nu }~\subset ~{\rm Range}\left( n v_{n,\ux }(\xi )\Big|_{|\xi |<\delta}~\right),   
$$   
as soon as $n$ is sufficiently large (also uniformly in $\ux\in\Space$), we are   
entitled to introduce the notation   
\begin{equation}   
\label{tilt_choice}   
\xi_{n,\ux } = \xi_{n,\ux }(r)~=~v_{n,\ux }^{-1}(r/n )\qquad\text{or, equivalently,}   
\qquad \frac{r}{n}~=~\nabla\Hnx (\xi_{n,\ux }) .   
\end{equation}   
Then the analytic implicit function theorem (c.f.~\cite{DS}) implies that, uniformly   
in $r\in R_{n,\nu }$ and $\ux\in\Space$,   
\begin{equation}   
\label{xi_n_x}   
\xi_{n,\ux }(r)~=~A_{\bf S}^{-1}\lb \frac{r}n-\nabla\log\rho_{\bf S}(0)\rb ~+~   
{\rm O}\lb n^{-2\nu}\rb .   
\end{equation}   
As a result, we conclude that, uniformly in $\ux\in\Space$ and $r\in R_{n,\nu}$,   
\begin{equation*}   
\begin{split}   
&\Qnx (r)~=~{\rm exp}\lbr -n\lb (\frac rn,\xi_{n,\ux}(r))_d -\Hnx (\xi_{n,\ux})\rb\rbr   
\bbQ_{n,\ux}^{\xi_{n, \ux}} (r)\\   
&\stackrel{\eqref{H_n_expansion},\eqref{xi_n_x}}{=}~   
{\rm exp}\lbr -\frac1{2n}\lb A_{\bf S}^{-1}(r-n\nabla\log\rho_{\bf S}(0)),  
(r-n\nabla\log\rho_{\bf S}(0))\rb_d\rbr  
\bbQ_{n,\ux}^{\xi_{n,\ux}} (r)\lb 1+\so \rb .   
\end{split}   
\end{equation*}   
Finally, by the very choice of the tilt $\xi_{n,\ux }(r)$ in \eqref{tilt_choice},   
the results of {\bf Step~1} apply to yield the desirable prefactor expression  
for $\bbQ_{n,\ux}^{\xi_{n,\ux}} (r)$.   
\qed

\subsection{Decay off the real axis}   
\label{sub_decay_away}   
In this subsection we establish the claim of Proposition~\ref{pro:off_gap}.   
 The proof involves   
three steps:   
\vskip 0.1cm  
\noindent   
\step{1} Fredholm spectrum of $L_{i\tau}$   
   
This has been already performed in Subsection~\ref{sub_gap}, and by   
\eqref{spectral_gap} of Theorem~\ref{thm:perron}, $\rho_{\bf F}\lb   
L_{i\tau}\rb~<~1$ holds for every $\tau\in [-\pi ,\pi ]^d$.   
\vskip 0.1cm   
\noindent   
\step{2} Spectrum of $L_{i\tau}$ for $\tau\neq 0$.   
\begin{lem}   
\label{lem:tau_gap}   
Assume that $\tau\neq 0$. Then,   
\begin{equation}   
\label{tau_gap}   
\rho_{\bf S}\lb L_{i\tau}\rb  ~<~1 .   
\end{equation}   
\end{lem}   
\noindent   
{\em Proof.} If there exists $\lambda\in\Sigma_{\bf S}\lb L_{i\tau }\rb$ with $\abs{\lambda}\geq 1$,   
then, by the preceding step, $\lambda$ is a Fredholm point and, as such,   
is, necessarily, an eigenvalue of $L_{i\tau}$. Let $h_\lambda\in\Ggt$ be a  
corresponding   
eigenfunction;   
$$   
L_{i\tau}h_\lambda ~=~\lambda h_\lambda .   
$$   
Taking the absolute values,   
$$   
L\abs{h_\lambda }~\geq ~\abs{\lambda}\,\abs{h_\lambda} .   
$$   
Since $L$ is normalized we, following the line of reasoning employed in  
Subsection~\ref{sub_pole}, infer that $\abs{\lambda}=1$ as well as   
that $h_\lambda$ is unimodal, $\abs{h_\lambda}\equiv 1$. Consequently, for every   
$\ux\in\Space$, every $n\in\bbN$ and each $\uz\in\cS_n$,   
\begin{equation}   
\label{z_indep}   
{\rm e}^{i\sum_1^n (\tau ,V(z_k ))_d}h_\lambda (\uz ,\ux )~=~\lambda^n h_\lambda (\ux ),   
\end{equation}   
and then, taking $n\to\infty$, conclude that $h_\lambda$ is a multiple of $\1$.   
In view of \eqref{z_indep} this means that $(\tau ,V(z))_d$ is independent of $z\in S$,  
a contradiction to the Assumption~{\bf A1} of Subsection~\ref{sub_local_setup}.   
\qed   
\vskip 0.1cm   
\noindent   
\step{3} Uniform estimate on $\rho_{\bf S}\lb L_{i\tau }\rb$.   
   
It remain to show that, given $\delta >0$, the inequality  
 \eqref{tau_gap} holds uniformly over  
$\tau\in [-\pi ,\pi]^d\setminus (-\gd ,\gd)^d$.This follows from   well known  
facts on the  
lower-semicontinuity of the spectrum. Assume that this is not the case, and  
there exists a sequence  
$\{\tau_k\}\subset [-\pi ,\pi]^d\setminus (-\gd ,\gd)^d $   
  and a  
sequence of numbers $\lambda_k\in\Sigma_{\bf S}\lb L_{i\tau_k }\rb$, such that  
$$   
\lim_{k\to\infty } \abs{\lambda_k} ~=~1 .   
$$   
Without loss of the generality we may assume that $\{\tau_k\}$ converges   
to some $\tau\neq 0$ and $\{\lambda_k\}$ converges to some $\lambda$ with  
$\abs{\lambda}=1$.   
By Lemma~\ref{lem:tau_gap}, however, $\rho_{\bf S}\lb L_{i\tau }\rb <1$. Therefore,  
$\lambda$   
belongs to the resolvent set of $ L_{i\tau }$. The latter is open, and one  
can find an $\epsilon >0$, such that the operator norm   
\begin{equation}   
\label{L_epsilon}   
\normgt{\lb \mu I-L_{i\tau}\rb^{-1}}~\leq ~\epsilon^{-1}   
\end{equation}   
for every $\abs{\mu -\lambda}\leq \epsilon$. On the other hand, $L_{i\tau_k}$ converges   
 to $L_{i\tau}$ in the strong operator topology: For every $f\in\Ggt$,   
$$   
\normsup{\lb L_{i\tau_k} - L_{i\tau}\rb f} ~\leq ~\normsup{f}\sum_{z\in S}   
\bigabs{1-{\rm e}^{i(\tau -\tau_k ,V(z))_d}} \phi (z),   
$$   
where we have introduced the notation   
$$   
\phi (z)~=~\sup_{\ux}{\rm e}^{\psi (z ,\ux)}\qquad\lb\text{Notice that by}   
\  \eqref{Psin_norm},\ \ \ \sum_z\phi (z) <\infty \rb .   
$$   
Similarly, using \eqref{Psin_norm},  
$$   
\snormgt{\lb L_{i\tau_k} - L_{i\tau}\rb f} ~\leq ~   
\lb c_3(\psi )\normsup{f} +\snormgt{f}\rb \sum_{z\in S}   
\bigabs{1-{\rm e}^{i(\tau -\tau_k ,V(z))_d}} \phi (z),   
$$   
with   
$$   
c_3(\psi )~=~\sup_{t\in (0,1)}\frac {{\rm e}^{\beta t} -1}{t} ,   
$$   
and $\beta$ specified in \eqref{beta_psi}. Thus,  
$\lim_{k\to\infty}\normgt{L_{i\tau_k }-L_{i\tau}}=0$,  
follows by the bounded convergence theorem.   
As a result, it follows from \eqref{L_epsilon} that $\lb  \mu I -L_{i\tau_k}\rb $ is  
invertible on $\setof{\mu}{\abs{\mu -\lambda}<\epsilon}$, as soon as $\tau_k$ is close  
enough to $\tau$, which is, of course, a  
contradiction.   
\qed   
   
\subsection{Perturbation theory of non-degenerate eigenvalues}   
\label{sub_perturbation}   
Let $\cF$ be a Banach space, $\cF^*$ its dual  
 and $B_\gd\subset \bbC$ is an open ball $B_\gd =   
\setof{z}{\abs{z}<\gd}$.

\noindent   
{\bf Definition.} A uniformly bounded  
family of linear operators $\{ T(\xi )\}_{\xi\in D}$ is   
said to be holomorphic on $B_\gd$ if   
$$   
\forall f\in\cF\ \text{and}\ \forall f^*\in\cF^*\qquad\text{the map}\ \xi\mapsto   
\lb T(\xi )f ,f^*\rb\ \text{is holomorphic in}\ B_\gd .  
$$   
We rely on the following statement of the analytic perturbation theory  
(c.f. \cite{Ka}, Section~VII.1.3 ):   
   
\vskip 0.1cm   
\noindent   
Let $\{ T(\xi )\}$ be a holomorphic family of operators on $B_\gd$, and assume that   
$\gl =\gl (0)$ is a  
{\em  non-degenerate} eigenvalue of $T(0)$. Then given a closed  
contour $\Gamma$ with ${\rm ext}(\Gamma )\cap\Sigma_{\bf S}\lb T(0)\rb =\{\gl\}$,   
there exists $\epsilon\in (0,\gd )$, such that:

\noindent   
1) For every $T(\xi )$ with  $\abs{\xi }<\epsilon$, there is exactly one  
spectral point $\gl (\xi )$, such that  
$\{ \gl (\xi )\}\, =\, {\rm ext}(\Gamma )\cap\Sigma_{\bf S}\lb T(\xi )\rb $.

\noindent   
2) $\gl (\xi )$ is a non-degenerate eigenvalue of $T(\xi )$ and the map  
$\xi\mapsto \gl (\xi )$ is analytic on $B_\epsilon$.   
\skip 0.1in   
\noindent   
We use this result in the following way:  By \eqref{Cramer_x} the family of  
operators $\{ L_\xi \}$ on $\Ggt$ is holomorphic  on $B_\gd$ for some $\gd >0$.   
According to  Theorem~\ref{thm:perron}, $\gl (0)=\rho_{\bf S}(L)$ is a non-degenerate   
eigenvalue of $L=L_0$ and, moreover, there exists $\nu >0$, such that the  
exteriour of $\Gamma_\nu \df\setof{z\in\bbC}{\abs{z}=(1-2\nu)\rho_{\bf S}(L)}$ satisfies   
$$   
\text{ext}\lb \Gamma_\nu\rb\cap\Sigma_{\bf S}\lb L_0 \rb~=~\{\gl (0)\} .   
$$   
Consequently, there exists $\epsilon >0$ and an analytic function  
$\lambda (\xi )$ on $B_\epsilon$, such that for every $\abs{\xi}\leq \epsilon$ the  
number $\gl (\xi )\df \rho_{\bf S}(L_\xi )$ is a non-degenerate eigenvalue of $L_\xi$,   
 $\abs{\gl (\xi )-\gl (0)}<\nu $, and  
$$   
\text{ext}\lb \Gamma_\nu\rb\cap\Sigma_{\bf S}\lb L_\xi \rb~=~\{\gl (\xi )\} .   
$$   
It follows that the family of projectors   
$$   
P_{L_\xi}~\df ~I~ +~\frac1{2\pi i}\oint_{\Gamma_\nu}\lb \mu I - L_\xi\rb^{-1}{\rm d}\mu   
$$   
is analytic on $B_\epsilon$, and so is the family  
$$   
\chi_{\ux} (\xi )\df P_{L_\xi}g (\ux ) ,   
$$   
 which shows up in the statement of Lemma~\ref{lem:perturbation}. Since  
$\inf_{\ux}   
 L g(\ux ) >0$, we can, if necessary, choose $\epsilon$ so small that   
$\{ \chi_{\ux }(\xi )\}$ is, uniformly in $\ux\in\Space$ and $\xi\in B_\epsilon$,   
 bounded away from zero. All the conclusions of Lemma~\ref{lem:perturbation} are,  
thereby, verified. Furthermore, for every $\xi\in B_\epsilon$ and $\ux\in\Space$,   
$$   
L^n_\xi  g (\ux)~=~ \chi_{\ux }(\xi )\rho_{\bf S}^n(\xi)~+~\smallo\lb (1-\nu)^n   
\abs{\rho_{\bf S}(\xi)}^n\rb .   
$$   
The expansion \eqref{H_n_expansion} follows.

\subsection{Non-degeneracy of ${\rm Hess}\lb \log\rho_{\bf S} \rb (0)$.}   
\label{sub_hessian}   
One has to show that there exists a positive $\alpha >0$, such that the variance  
$$   
\lb {\rm Hess}(\Hnx )l,l\rb_d~=~\frac1n{\bbV}{\rm ar}_{n,\ux}\lb  
\sum_{k=1}^n \lb V(z_i ),l\rb_d\rb~\geq ~\alpha \abs{l}^2,   
$$   
uniformly in $n$, $\ux\in\Space$ and $l\in\bbR^d$. This follows from the conditional   
variance argument based on the assumptions {\bf A1},{\bf A2} and the   
H\"{o}lder upper bound  \eqref{Psin_norm}. Indeed, let $n=km+l$. Then,   
$$   
\bbV{\rm ar}_{n,\ux}\lb \sum_{i=1}^n \lb V(z_i ),l\rb_d\rb~   
\geq ~\min_{\bar{z}_j\in S\, j\neq 0\,{\rm mod}(m )}   
\bbV{\rm ar}_{n,\ux}\lb \sum_{i=1}^{k} \lb V(z_{im}),l\rb_d \sep   
z_j =\bar{z}_j\,  
\rb   
$$   
However the conditional variances of the variables $V(z_{im} )$ are,  
uniformly in $\{ \bar{z}_j\}$, bounded away both from zero and $\infty$, whereas  
the correlation coefficient between different $ V(z_{im} )$'s decays to zero   
exponentially fast with $m$.\qed

\end{document}